\documentclass[11pt]{article}
\usepackage{times,mathptm}
\usepackage{amsfonts,amscd,amssymb}
\usepackage{theorem}
\input{diagrams}
\newtheorem{definition}{Definition}[section]
\newtheorem{lemma}[definition]{Lemma}
\newtheorem{proposition}[definition]{Proposition}
\newtheorem{corollary}[definition]{Corollary}
{\theorembodyfont{\rmfamily}\newtheorem{remark}[definition]{Remark}}
{\theorembodyfont{\rmfamily}\newtheorem{remarks}[definition]{Remarks}}
\newtheorem{theorem}[definition]{Theorem}
{\theorembodyfont{\rmfamily}}
{\theorembodyfont{\rmfamily}}
{\theorembodyfont{\rmfamily}}
\def\va{\varepsilon}
\def\v{\varphi}

\def\tr{\triangleright}
\def\rh{\rightharpoonup}
\def\lh{\leftharpoonup}

\def\ra{\rightarrow}
\def\a{\alpha}
\def\b{\beta}

\def\l{\lambda}
\def\r{\rho}
\def\cd{\cdot}
\def\d{\delta}
\def\ov{\overline}
\def\un{\underline}
\newcommand{\ovc}{\mbox{$M^{\ov {{\rm co}H}}$}}
\newcommand{\ove}{\mbox{$\ov {E}$}}
\newcommand{\hb}{\mbox{$_H{\cal M}_H^H$}}
\newcommand{\smi}{\mbox{$S^{-1}$}}

\def\rawo\lonra{\longrightarrow}

\def\ot{\otimes}
\def\act{\hspace*{-1pt}\ra\hspace*{-1pt}}

\newcommand{\thlabel}[1]{\label{th:#1}}
\newcommand{\thref}[1]{Theorem~\ref{th:#1}}
\newcommand{\selabel}[1]{\label{se:#1}}
\newcommand{\seref}[1]{Section~\ref{se:#1}}
\newcommand{\lelabel}[1]{\label{le:#1}}
\newcommand{\leref}[1]{Lemma~\ref{le:#1}}
\newcommand{\prlabel}[1]{\label{pr:#1}}
\newcommand{\prref}[1]{Proposition~\ref{pr:#1}}
\newcommand{\colabel}[1]{\label{co:#1}}
\newcommand{\coref}[1]{Corollary~\ref{co:#1}}
\newcommand{\relabel}[1]{\label{re:#1}}

\newcommand{\delabel}[1]{\label{de:#1}}

\newcommand{\eqref}[1]{(\ref{eq:#1})}

\newenvironment{proof}{{\it Proof.}}{\hfill $ \square $ \vskip 4mm}
\begin{document}
\title{Integrals for (dual) quasi-Hopf algebras. Applications
\thanks{Research supported by the bilateral project
``Hopf Algebras in Algebra, Topology, Geometry and Physics" of the Flemish and
Romanian governments.}}
\author{D. Bulacu\thanks{This paper was written while the first author was
visiting
the Vrije Universiteit Brussel, and he would like to
thank VUB for its warm hospitality.}\\
Faculty of Mathematics\\ University of Bucharest\\
RO-70109 Bucharest 1, Romania\and
S. Caenepeel
\\ Faculty of Applied Sciences\\
Vrije Universiteit Brussel, VUB\\ B-1050 Brussels, Belgium}
\date{}
\maketitle

\begin{abstract}
A classical result in the theory of Hopf algebras concerns the
uniqueness and existence of integrals: for an arbitrary Hopf algebra,
the integral space has dimension $\leq 1$, and for a finite
dimensional Hopf algebra, this dimension is exaclty one. We
generalize these results to quasi-Hopf algebras and dual
quasi-Hopf algebras. In particular, it will follow that the
bijectivity of the antipode follows from the other axioms of
a finite dimensional quasi-Hopf algebra. We give a new version
of the Fundamental Theorem for quasi-Hopf algebras. We show that
a dual quasi-Hopf algebra is co-Frobenius if and only if it has
a non-zero integral. In this case, the space of left or right
integrals has dimension one.
\end{abstract}

\section{Introduction}\selabel{0}
Quasi-bialgebras and quasi-Hopf algebras were introduced by
Drinfeld in \cite{d1}, in connection with the Knizhnik-Zamolodchikov
system of partial differential equations, cf. \cite{Kassel95}.
From a categorical point of view, the
notion is not so different from classical bialgebras: we consider
an algebra $H$, and we want to make the category of $H$-modules,
equipped with the tensor product of vector spaces, into a
monoidal category. If we require that the associativity constraint
is the natural associativity condition for vector spaces, then we
obtain a bialgebra structure on $H$, in general, we obtain a
quasi-bialgebra structure, that is, we have a comultiplication
and a counit on $H$, where the comultiplication is not
necessarily coassociative, but only quasi-coassociative.\\
Of course the theory of quasi-bialgebras and quasi-Hopf algebras
is technically more complicated than the classical Hopf algebra
theory. A more conceptual difference however, is the fact that
the definition of a bialgebra is self-dual, and this symmetry
is broken when we pass to quasi-bialgebras. As a consequence,
we don't have the notion of comodule or Hopf module over a
quasi-Hopf algebra, and results in Hopf algebras that depend
on these notions cannot be generalized in a straightforward way.
For instance, the classical proof of the uniqueness and
existence of integral is based on the Fundamental Theorem
for Hopf modules \cite{sw}.\\
Hausser and Nill \cite{hn3} proved that a finite dimensional quasi-Hopf
algebra is a Frobenius algebra, and has a one dimensional integral
space. Independently, Panaite and Van Oystaeyen \cite{pv2} proved
the existence of integrals for finite dimensional quasi-Hopf algebras,
using the approach developed in \cite{vd}, without using
quasi-Hopf bimodules.\\
For a finite dimensional Hopf algebra $H$, it follows from the Fundamental
Theorem that $\int _l^H\ot H^*$ and $H$ are isomorphic as Hopf modules.
In \seref{2}, we will see that the isomorphism survives as a
left $H$-linear isomorphism in the case of a finite dimensional
quasi-Hopf algebra. The method of proof is quite different from the
classical one: the isomorphism is constructed explicitly, using the
projection of $H$ onto the integral space constructed in \cite{pv2}.
In Drinfeld's original definition \cite{d1}, the antipode of a quasi-Hopf
algebra is required to be bijective. Actually our proof of
 $\int _l^H\ot H^* \cong H$ does not use this bijectivity, and has as
a consequence that, for a finite dimensional quasi-Hopf algebra,
the bijectivity of the antipode follows from the other axioms; another
consequence is that the integral space is one-dimensional. In
a recent preprint \cite{sh}, Schauenburg gave
a different proof of the fact that the antipode of a finite
dimensional quasi-Hopf Hopf algebra is bijective. \\
The infinite dimensional case is treated as well. We show that a quasi-Hopf
algebra (without the assumption that the antipode is bijective) is
finite dimensional if and only if the antipode is bijective and the
integral space is non-zero. The integral space of an infinite dimensional
quasi-Hopf algebra with bijective antipode is zero. A semisimple Hopf
algebra with bijective antipode is finite dimensional.\\
Hausser and Nill \cite{hn3} also introduced cointegrals on a finite
dimensional quasi-Hopf algebra; these cointegrals are elements of the
dual space $H^*$, and, using a Structure Theorem for quasi-Hopf bimodules,
Hausser and Nill prove that the space of cointegrals ${\cal L}$ is
one-dimensional,
and that all non-zero integrals are nondegenerate. In \seref{3},
we further investigate cointegrals. In \cite{hn3}, it is asked whether
there is a connection between the projection of $H$ onto the space of
integrals from \cite{pv2}, and the projection of $H^*$ onto the space
of cointegrals, introduced in \cite{hn3}. This is done in \leref{3.2},
and, as an application, we give some characterizations of cointegrals,
see \prref{3.4}.\\
In the second part of \seref{3}, we propose an alternative definition of
the space of coinvariants of a quasi-Hopf bimodule. This alternative
space of coinvariants is isomorphic to the Hausser-Nill space of coinvariants,
and can be used to give a second version of the Structure Theorem. Our
alternative has nevertheless two advantages, compared to the Hausser-Nill
approach: first, it is invariant under the adjoint action (cf. \leref{3.5}).
Secondly, in the finite dimensional case, it gives rise to an alternative
definition of cointegral: we take the alternative coinvariants of $H^*$.
If we write down this formula explicitely, we obtain a formula that still
makes sense in the infinite dimensional case, so we obtain a plausible
definition for cointegrals in the infinite dimensional case.\\
As we have already pointed out, the definition of quasi-Hopf algebra is
not self-dual. Actually, we can introduce dual quasi-Hopf algebras, these
are coalgebras, with a multiplication that is not associative, but only
quasi-associative. In \seref{4}, we introduce integrals in
dual quasi-Hopf algebras. We were able to prove that the rational dual
of a dual quasi-Hopf algebra $A$ is isomorphic as a comodule
to the tensor product of $A$ itself and the integral space. This generalizes
the classical statement for Hopf algebras (see \cite{sw}), but, again,
we have to give a direct proof, and cannot deduce the statement from a
Structure Theorem. As in the classical case, it then follows immediately
that the integral space is zero if and only if the rational dual is zero.
Also we can use the integrals to investigate properties of a dual
quasi-Hopf algebra $A$ as a coalgebra (\thref{4.5}). The existence of a
non-zero integral is equivalent to $A$ being a co-Frobenius coalgebra,
a QcF coalgebra, or a left semiperfect coalgebra. Moreover, for a
dual quasi-Hopf algebra, all these notions are left-right symmetric.
Furthermore the existence of a non-zero integral is equivalent to $A$
being a generator or a projective object in the category of (left or
right) comodules.\\
As a first application of this coalgebraic viewpoint, we find that
a dual quasi-Hopf subalgebra of a dual quasi-Hopf algebra with non-zero
integral has non-zero integrals. Secondly, it follows that non-zero
integrals are unique up to multiplication by a scalar. Also we can
give the connection between left and right integrals (\prref{4.9}),
and this generalizes \cite[Proposition 1.3]{bdgn}. We were able to prove
that the antipode of a dual quasi-Hopf algebra with a non-zero integral is
injective,
but it remains open if it is also surjective, as it is the case for
a classical Hopf algebra, see \cite{re}. Our final result is Maschke's
Theorem for dual
quasi-Hopf algebras (\thref{4.10}), stating that a dual quasi-Hopf algebra is
cosemisimple if and only there exists an integral $T$ such that $T(1)=1$.\\
When we pass from bialgebras and Hopf algebras to
quasi-bialgebras and quasi-Hopf algebras, the appearance of the
reassociator and the elements $\alpha$ and $\beta$ in the definition
of the antipode, considerably increase the complexity of computations
and proofs. This observation is not new, other authors who have been
working on quasi-Hopf algebras experienced this before us. However, the
philosophy is basically the same as in the case of usual bialgebras: the
idea is to make the category of $A$-modules into a monoidal category.
Recently, Schauenburg proposed an alternative approach to proving
results on quasi-bialgebras, exploiting the categorical ideas behind
quasi-bialgebras, and replacing the computational arguments using the
Sweedler notation by conceptual arguments (see \cite{sh} for detail). At
this moment is
not clear to us whether Schauenburg's ideas can be used to give alternative
and/or
more transparant proofs of the results in this paper.

\section{Preliminaries}\selabel{1}
We work over a commutative field $k$. All algebras, linear
spaces etc. will be over $k$; unadorned $\ot $ means $\ot_k$.
Following Drinfeld \cite{d1}, a quasi-bialgebra is
a four-tuple $(H, \Delta
, \va , \Phi )$ where $H$ is an associative algebra with unit,
$\Phi$ is an invertible element in $H\ot H\ot H$, and
$\Delta :\
H\ra H\ot H$ and $\va :\ H\ra k$ are algebra homomorphisms satisfying
the identities
\begin{equation}\label{q1}
(id \ot \Delta )(\Delta (h))=
\Phi (\Delta \ot id)(\Delta (h))\Phi ^{-1},
\end{equation}
\begin{equation}\label{q2}
(id \ot \va )(\Delta (h))=h\ot 1,
\mbox{${\;\;\;}$}
(\va \ot id)(\Delta (h))=1\ot h,
\end{equation}
for all $h\in H$, and
$\Phi$ has to be a $3$-cocycle, in the sense that
\begin{equation}\label{q3}
(1\ot \Phi)(id\ot \Delta \ot id)
(\Phi)(\Phi \ot 1)=
(id\ot id \ot \Delta )(\Phi )
(\Delta \ot id \ot
id)(\Phi ),
\end{equation}
\begin{equation}\label{q4}
(id \ot \va \ot id
)(\Phi )=1\ot 1\ot 1.
\end{equation}
The map $\Delta $ is called the
coproduct or the
comultiplication, $\va $ the counit and $\Phi $ the
reassociator.
As for Hopf algebras we denote $\Delta (h)=\sum h_1\ot h_2$,
but
since $\Delta$ is only quasi-coassociative we adopt the
further
convention
$$ 
(\Delta \ot id)(\Delta (h))= \sum h_{(1,
1)}\ot
h_{(1, 2)}\ot h_2, \mbox{${\;\;\;}$} (id\ot \Delta )(\Delta
(h))=
\sum h_1\ot h_{(2, 1)}\ot h_{(2,2)}, 
$$
for all $h\in H$. We will
denote the tensor components of $\Phi$
by capital letters, and the ones of
$\Phi^{-1}$ by small letters,
namely
\begin{eqnarray*}
&&\Phi=\sum X^1\ot
X^2\ot X^3= \sum T^1\ot T^2\ot T^3=
\sum V^1\ot V^2\ot
V^3=\cdots\\
&&\Phi^{-1}=\sum
x^1\ot x^2\ot x^3= \sum t^1\ot t^2\ot t^3=
\sum v^1\ot v^2\ot v^3=\cdots
\end{eqnarray*}
$H$ is
called a quasi-Hopf
algebra if, moreover, there exists an
anti-automorphism $S$ of the algebra
$H$ and elements $\a , \b \in
H$ such that, for all $h\in H$, we
have:
\begin{equation}\label{q5}
\sum S(h_1)\a h_2=\va
(h)\a
\mbox{${\;\;\;}$ and ${\;\;\;}$}
\sum h_1\b S(h_2)=\va (h)\b
,
\end{equation}
\begin{equation} \label{q6}
\sum X^1\b S(X^2)\a
X^3=1
\mbox{${\;\;\;}$ and${\;\;\;}$}
\sum S(x^1)\a x^2\b
S(x^3)=1.
\end{equation}
For a quasi-Hopf algebra the antipode is
determined uniquely up
to a transformation $\a \mapsto U\a $, $\b \mapsto
\b U^{-1}$,
$S(h)\mapsto US(h)U^{-1}$, where $U\in H$ is invertible.
The
axioms for a quasi-Hopf algebra imply that $\va (\a )\va (\b )=1$,
so, by
rescaling $\a $ and $\b $, we may assume without loss of generality
that
$\va (\a )=\va (\b )=1$ and $\va \circ S=\va $. The identities
(\ref{q2}),
(\ref{q3}) and (\ref{q4}) also imply that
\begin{equation}\label{q7}
(\va
\ot id\ot id)(\Phi )=
(id \ot id\ot \va )(\Phi )=1\ot 1\ot
1.
\end{equation}
Next we recall that the definition of a
quasi-Hopf
algebra is ``twist coinvariant" in the following sense.
An
invertible element $F\in H\ot H$ is called a {\sl gauge
transformation}
or {\sl twist} if
$(\va \ot id)(F)=(id\ot \va)(F)=1$.
If
$H$ is a quasi-Hopf algebra and $F=\sum F^1\ot F^2\in H\ot H$
is a gauge
transformation with inverse $F^{-1}=\sum G^1\ot
G^2$, then we can define a
new quasi-Hopf algebra $H_F$ by keeping
the multiplication, unit, counit
and antipode of $H$ and
replacing the comultiplication, antipode and the
elements $\alpha$
and $\beta$ by
\begin{equation} \label{g1}
\Delta
_F(h)=F\Delta (h)F^{-1},
\end{equation}
\begin{equation}
\label{g2}
\Phi_F=(1\ot F)(id \ot \Delta )(F) \Phi (\Delta \ot
id)
(F^{-1})(F^{-1}\ot 1),
\end{equation}
\begin{equation}
\label{g3}
\a_F=\sum S(G^1)\a G^2,
\mbox{${\;\;\;}$}%
\b_F=\sum F^1\b
S(F^2).
\end{equation}
It is well-known that the antipode of a Hopf algebra
is an
anti-coalgebra morphism. For a quasi-Hopf algebra, we have
the
following statement: there exists a gauge transformation
$f\in H\ot H$ such
that
\begin{equation} \label{ca}
f\Delta (S(h))f^{-1}= \sum (S\ot S)(\Delta
^{op}(h))
\mbox{,${\;\;\;}$for all $h\in H$.}
\end{equation}
$f$ can be
computed explicitly. First set
\begin{equation}
\sum A^1\ot A^2\ot A^3\ot
A^4=
(\Phi \ot 1) (\Delta \ot id\ot id)(\Phi
^{-1}),
\end{equation}
\begin{equation} \sum B^1\ot B^2\ot B^3\ot
B^4=
(\Delta \ot id\ot id)(\Phi )(\Phi ^{-1}\ot 1)
\end{equation}
and then
define $\gamma, \delta\in H\ot H$by
\begin{equation} \label{gd}%
\gamma
=\sum S(A^2)\a A^3\ot S(A^1)\a A^4~~{\rm and}~~
\delta
=\sum B^1\b
S(B^4)\ot B^2\b S(B^3).
\end{equation}
$f$ and $f^{-1}$ are then given by
the formulas
\begin{eqnarray}
f&=&\sum (S\ot S)(\Delta ^{op}(x^1)) \gamma
\Delta (x^2\b
S(x^3)),\label{f}\\
f^{-1}&=&\sum \Delta (S(x^1)\a x^2)
\delta
(S\ot S)(\Delta ^{op}(x^3))\label{g}
\end{eqnarray}
where $\Delta
^{op}(h)=\sum h_2\ot h_1$.
$f$ satisfies the following
relations:
\begin{equation} \label{gdf}
f\Delta (\a )=\gamma
,
\mbox{${\;\;\;}$}
\Delta (\b )f^{-1}=\delta .
\end{equation}
Furthermore
the corresponding
twisted reassociator (see (\ref{g2})) is given
by
\begin{equation} \label{pf}
\Phi _f=\sum (S\ot S\ot S)(X^3\ot X^2\ot
X^1).
\end{equation}
In a Hopf algebra $H$, we obviously have the
identity
$$
\sum h_1\ot h_2S(h_3)=h\ot 1,~{\rm for~all~}h\in H.
$$
We will
need the generalization of this formula to
quasi-Hopf algebras. Following
\cite{hn1}, \cite{hn2}, we define
\begin{equation} \label{qr}
p_R=\sum
p^1\ot p^2=\sum x^1\ot x^2\b S(x^3),
\mbox{${\;\;\;}$} q_R=\sum q^1\ot
q^2=\sum X^1\ot S^{-1}(\a X^3)X^2.
\end{equation}
For all $h\in H$, we then
have
\begin{equation} \label{qr1}
\sum \Delta (h_1)p_R[1\ot
S(h_2)]=p_R[h\ot 1],
\mbox{${\;\;\;}$}%
\sum [1\ot S^{-1}(h_2)]q_R\Delta
(h_1)=(h\ot 1)q_R,
\end{equation}
and
\begin{equation} \label{pqr}
\sum
\Delta (q^1)p_R[1\ot S(q^2)]=1\ot 1,
\mbox{${\;\;\;}$}
\sum [1\ot
S^{-1}(p^2)]q_R\Delta (p^1)=1\ot
1,
\end{equation}
\begin{eqnarray}
&&\hspace*{-2cm}(q_R\ot 1)(\Delta \ot
id)(q_R)\Phi ^{-1}\nonumber\\
&=&\sum [1\ot S^{-1}(X^3)\ot S^{-1}(X^2)]
[1\ot S^{-1}(f^2)\ot
S^{-1}(f^1)] (id \ot \Delta )(q_R\Delta
(X^1)),\label{qr2}\\
&&\hspace*{-2cm}\Phi (\Delta \ot id)(p_R)(p_R\ot
id)\nonumber\\
&=&\sum (id\ot \Delta )(\Delta (x^1)p_R)(1\ot f^{-1})(1\ot
S(x^3)\ot
S(x^2))\label{pr},
\end{eqnarray}
where $f=\sum f^1\ot f^2$ is
the twist defined in (\ref{f}). Note
that the formulas (\ref{g}-\ref{pr})
(except (\ref{qr2}) and
the second part of (\ref{qr1}-\ref{pqr})) can be
proved
without using the bijectivity of the antipode
$S$.

\section{Integrals in quasi-Hopf
algebras}\selabel{2}
\setcounter{equation}{0}
Let $H$ be a finite
dimensional quasi-Hopf algebra with an
antipode $S$. In \cite{hn3}, it is
shown that
$H$ is a Frobenius algebra and, as a consequence,
the space of
left (right) integrals in $H$ is one dimensional.
The proof in \cite{hn3}
relies on the fact that the antipode
$S$ is bijective. Using different
arguments independent of
the bijectivity of the antipode, we will give
another
proof of the existence and uniqueness of integrals in $H$.
In fact
we will prove that, in the definition of finite
dimensional quasi-Hopf
algebra, the bijectivity of the antipode
can be dropped, in other words,
the bijectivity of the
antipode follows from the other axioms. This will
generalize
a similar result for Hopf algebras, see \cite{ls}.\\
Let us make
our terminology consistent: by a quasi-Hopf algebra,
we will mean a
quasi-Hopf algebra as defined in \seref{1},
but without the assumption that
the antipode is bijective. If
the antipode is bijective, then we will say
that we have
a quasi-Hopf algebra in the sense of Drinfeld. With this
convention,
our main result is the following: a finite dimensional
quasi-Hopf
algebra is a quasi-Hopf algebra in the sense of
Drinfeld.\\
Recall
that $t\in H$ is called a left (resp. right) integral in
$H$
if $ht=\va (h)t$ (resp. $th=\va (h)t)$, $\forall h\in H$.
We denote by
$\int _l^H$ ($\int _r^H$) the space of left (right)
integrals in $H$. If
there exists a non-zero left integral in $H$ which is
at the same time a
right integral, then $H$ is called
unimodular.\\
In \cite{ls}, the
Fundamental Theorem is used to prove the existence
and uniqueness of
integrals, and then the bijectivity of the
antipode follows. In the
quasi-Hopf algebra case, this approach
will not work, since we cannot
define Hopf modules. Van Daele
\cite{vd} gave a short and direct proof of
the existence and
uniqueness of integrals in a finite dimensional Hopf
algebra,
and Panaite and Van Oystaeyen
\cite{pv2} generalized Van Daele's
argument, proving the existence
of left integrals in finite dimensional
quasi-Hopf algebras.
More precisely,
let $H$ be a finite dimensional
quasi-Hopf algebra,
$\{e_i\}_{i=\ov {1, n}}$ a basis of $H$ and
$\{e^i\}_{i=\ov {1,
n}}$ the dual basis of $H^*$. Following \cite{pv2}
we
define
\begin{equation}\label{p}
P(h)=\sum \limits _{i=1}^n\langle e^i, \b
S(S(X^2(e_i)_2)\a X^3)h\rangle
X^1(e_i)_1,
\end{equation}
for all $h\in H$.
Then one can show that
$P(h)\in \int _l^H, \forall h\in H$ and $\sum
\limits _{i=1}^n
\langle e^i, S(P(e_i)\b )\rangle =\va (\b )=1$. It follows
that at least
one of the
$P(e_i)\neq 0$, and $\int_l^H\not=0$.\\
In order
to prove the uniqueness of integrals for finite
dimensional quasi-Hopf
algebras we need the following

\begin{lemma}\lelabel{2.1}
Let $t$ be a
left integral in a quasi-Hopf algebra $H$. Then
for all $h\in
H$
\begin{equation}\label{f1}
\sum hX^1t_1\ot S(X^2t_2)\a X^3=\sum
X^1t_1\ot S(X^2t_2)\a
X^3h,
\end{equation}
and
\begin{equation}\label{f2}
\sum t_1\ot S(t_2)=\sum
X^1t_1\ot  S(X^2t_2)\a X^3\b =
\sum \b X^1t_1\ot S(X^2t_2)\a X^3
.
\end{equation}
\end{lemma}

\begin{proof}
For all $h\in H$ we calculate,
using (\ref{q5}), (\ref{q1}) and
$t\in \int _l^H$
that
\begin{eqnarray*}
&&\hspace*{-2cm}\sum hX^1t_1\ot S(X^2t_2)\a
X^3=
\sum h_1X^1t_1\ot S(h_{(2, 1)}X^2t_2)\a h_{(2, 2)}X^3
\\
&=&\sum
X^1(h_1t)_1\ot S(X^2(h_1t)_2)\a X^3h_2
=\sum X^1t_1\ot S(X^2t_2)\a
X^3h.
\end{eqnarray*}
To prove the first equality in (\ref{f2}),
we take
$\sum X^1t_1\ot S(X^2t_2)\a X^3\b $. First we apply
the 3-cocycle
condition
$$
\Phi \ot 1=(id\ot \Delta \ot id)(\Phi ^{-1})(1\ot \Phi
^{-1})
(id\ot id\ot \Delta )(\Phi )(\Delta \ot id\ot id)(\Phi )
$$
and
then, successively using the fact that
$t\in \int _l^H$,
(\ref{q5}),
(\ref{q4}), (\ref{q7}) and (\ref{q6}),
we find the left hand side of
(\ref{f2}). The second equality
in (\ref{f2}) follows from
(\ref{f1}).
\end{proof}

For a quasi-Hopf algebra $H$, we introduce $H^*$
as the dual space of
$H$ with its natural multiplication
$\langle h^*g^*,
h\rangle =\sum h^*(h_1)g^*(h_2)$, where $h^*, g^*\in H^*$
and $h\in H$. If $H$ is finite dimensional, then $H^*$
is also equipped with a natural coassociative
coalgebra structure $(\tilde {\Delta }, \tilde {\va })$
given by $\langle \tilde {\Delta }(h^*), h\ot h'\rangle =
\langle h^*, hh'\rangle
$ and $\tilde {\va }(h^*)=h^*(1)$, where
$h^*\in H^*$, $h, h'\in H$ and
$\langle  , \rangle :\ H^*\ot H\ra k$
denotes the dual pairing. On $H^*$ we
have the natural left and
right $H$-actions
\begin{equation}\label{act}
\langle h\rh h^*, h'\rangle
=\langle h^*, h'h\rangle,
\mbox{${\;\;}$}
\langle h^*\lh h, h'\rangle
=\langle h^*,hh'\rangle,
\end{equation}
where $h, h'\in H$ and $h^*\in
H^*$. This makes $H^*$
into a $H$-$H$-bimodule.\\
We also introduce
$\ov {S}:\ H^*\ra H^*$ as the anti-coalgebra homomorphism
dual to $S$, i.e.
$\langle \ov {S}(h^*), h\rangle =\langle h^*, S(h)\rangle $,
$\forall h^*\in H^*, h\in H$.

\begin{theorem}\thlabel{2.2}
Let $H$ be a finite
dimensional quasi-Hopf algebra, $\{e_i\}_{i=\ov {1, n}}$
a basis of $H$
with dual basis $\{e^i\}_{i=\ov {1, n}}$ of $H^*$, and
define $\theta :\
\int _l^H\ot H^*\ra H$, by
\begin{equation}\label{th}
\theta (t\ot
h^*)=\sum h^*(S(X^2t_2p^2)\a X^3)X^1t_1p^1,
\mbox{${\;\;}$$\forall t\in
\int _l^H, h^*\in H^*$,}
\end{equation}
where $p_R=\sum p^1\ot p^2$ is the
element defined in (\ref{qr}).
Then the following assertions hold:\\
(i)
$\theta $ is an isomorphism of left $H$-modules,
where $\int _l^H\ot H^*$
is a left $H$-modules via
$h\cd (t\ot h^*)=t\ot h\rh h^*$ $\forall h\in H,
t\in \int _l^H,
h^*\in H^*$, and $H$ is a left $H$-module via left
multiplication.
Consequently
$dim_k\int _l^H=1$. The inverse of $\theta$ is
given by
\begin{equation}\label{thi}
\theta^{-1}(h)=\sum \limits
_{i=1}^nP(e_ih)\ot e^i,
\mbox{${\;\;}$$\forall h\in
H$,}
\end{equation}
where $P$ is the projection onto the space of left
integrals, defined in
(\ref{p}).\\
(ii) The antipode $S$ is
bijective.\\
(iii) $S(\int _l^H)=\int _r^H$, $S(\int _r^H)=\int _l^H$,
and
$dim_k\int _r^H=1$.
\end{theorem}

\begin{proof} (i) First we show that
$\theta $ and $\theta ^{-1}$
are inverses. Indeed, for all $h\in H$ we
have:
\begin{eqnarray*}
\theta (\theta ^{-1}(h))&\pile{{\rm (\ref{th},
\ref{p})}\\ =}&
\sum \limits _{i, j=1}^n \langle e^j, \b
S(S(X^2(e_j)_2)\a
X^3)e_ih\rangle  \langle e^i, S(Y^2X^1_2(e_j)_{(1, 2)}p^2)
\a Y^3\rangle
\\
&&~~~~~Y^1X^1_1(e_j)_{(1, 1)}p^1 \\
{\rm (\ref{qr1})}~~~~&=&\sum \limits
_{j=1}^n\langle e^j, \b
S(Y^2X^1_2p^2S(X^2)\a
X^3)\a Y^3h\rangle
Y^1X^1_1p^1e_j\\
{\rm (\ref{qr}, \ref{q3})}~~~~&=&\sum Y^1x^1Z^1\b
S(Y^2x^2X^1Z^2_1\b
S(x^3_1X^2Z^2_2)\a x^3_2X^3Z^3) \a Y^3h\\
{\rm
(\ref{q5},\ref{q4},\ref{q7},\ref{q6})}~~~~&=&\sum Y^1\b S(Y^2X^1\b
S(X^2)\a
X^3)\a Y^3h=h.
\end{eqnarray*}
For all $t\in \int _l^H$ and $h^*\in H^*$ we
compute
\begin{eqnarray*}
\theta ^{-1}(\theta (t\ot h^*))
&=&\sum \limits
_{i=1}^nh^*(S(X^2t_2p^2)\a X^3)P(e_iX^1t_1p^1)\ot
e^i\\
{\rm (\ref{p},
\ref{f1})}~~~~&=&\sum \limits _{i, j=1}^n
h^*(S(X^2t_2p^2)\a X^3\b
S(S(Y^2(e_j)_2)\a Y^3)e_i)
\langle e^j, X^1t_1p^1\rangle Y^1(e_j)_1\ot
e^i\\
{\rm (\ref{f2})}~~~~&=&\sum \limits _{i=1}^n 
h^*(S(S(Y^2t_{(1,
2)}p^1_2)\a
Y^3t_2p^2)e_i)Y^1t_{(1, 1)}p^1_1\ot e^i\\
{\rm
(\ref{q1},\ref{q5})}~~~~&=&\sum \limits _{i=1}^n
h^*(S(S(Y^2p^1_2)\a
Y^3p^2)e_i)tY^1p^1_1\ot e^i\\
{\rm (\ref{qr}, \ref{q3})}~~~~&=&\sum \limits
_{i=1}^n 
h^*(S(S(y^2_1x^1X^2)\a
y^2_2x^2X^3_1\b
S(y^3x^3X^3_2))e_i)ty^1X^1\ot e^i\\
{\rm
(\ref{q5},\ref{q4},\ref{q7})}~~~~&=&\sum \limits _{i=1}^n
h^*(S(S(x^1)\a
x^2\b S(x^3))e_i)t\ot e^i\\
{\rm (\ref{q6})}~~~~&=&\sum \limits _{i=1}^n
h^*(e_i)t\ot e^i
=t\ot h^*.
\end{eqnarray*}
Since $\theta $ is a bijection
and ${\rm dim}_kH={\rm dim}_kH^*$ is finite,
is follows that ${\rm
dim}_k\int _l^H=1$. We are left to show that
$\theta $ is $H$-linear. For
all $h\in H$, $t\in \int
_l^H$ and $h^*\in H^*$ we
have:
\begin{eqnarray*}
h\theta (t\ot h^*)
&=&\sum h^*(S(X^2t_2p^2)\a
X^3)hX^1t_1p^1\\
{\rm (\ref{f1})}~~~~&=&\langle h\rh h^*, 
S(X^2t_2p^2)\a
X^3\rangle X^1t_1p^1\\
&=&\theta (t\ot h\rh h^*).
\end{eqnarray*}

(ii)
First we prove
that $\ov {S}$ is bijective. $H^*$ is finite dimensional, so
it suffices to show that $\ov{S}$ is injective. Let $h^*\in H^*$ be such
that $\ov {S}(h^*)=0$, and take $0\not=t\in \int _l^H$. For all $h\in H$ we
have
\begin{eqnarray*}
&&\hspace*{-1cm}\theta (t\ot \b S(h)\rh
h^*)
=\langle \b S(h)\rh h^*, S(X^2t_2p^2)\a X^3\rangle X^1t_1p^1\\
&=&\sum
\langle h^*, S(X^2t_2p^2)\a X^3\b S(h)\rangle X^1t_1p^1\\
{\rm
(\ref{f2})}~~~~&=&\sum \langle h^*, S(ht_2p^2)\rangle t_1p^1\\
&=&\sum
\langle \ov {S}(h^*), ht_2p^2\rangle t_1p^1=0.
\end{eqnarray*}
 Since
$\theta $ is bijective we obtain that
$t\ot \b S(h)\rh h^*=0$. Now, because
$t\not=0$ and
${\rm dim}_k\int _l^H=1$, it follows that
$\b S(h)\rh h^*=0$,
for all $h\in H$. Therefore, by (\ref{q6}), for
all $h'\in H$ we
have
$$
h^*(h')=\sum \langle h^*, h'S(x^1)\a x^2\b S(x^3)\rangle = \sum
\langle \b
S(x^3)\rh h^*, h'S(x^1)\a x^2\rangle =0.
$$
It is not hard to
see that
$\ov {S}^*:\ H^{**}\ra H^{**}$, $\ov {S}^*(h^{**})=
h^{**}\circ
\ov {S}$, $\forall h^{**}\in H^{**}$, is a bijective
map. If we define $\xi
:\ H\ra H^{**}$ by $\xi (h)(h^*)=h^*(h)$,
$\forall h\in H, h^*\in H^*$,
then we can easily show that $\{\xi
(e_i)\}_{i=\ov {1, n}}$ is a basis of
$H^{**}$ dual to the basis
$\{e^i\}_{i=\ov {1, n}}$ of $H^*$, and it
follows that $\xi $ is bijective.
Moreover, $\xi ^{-1}$ is given by $\xi
^{-1}(h^{**})=\sum \limits
_{i=1}^nh^{**}(e^i)e_i, \forall h^{**}\in
H^{**}$. In addition,
$\xi ^{-1}\circ \ov {S}^*\circ \xi =S$ so $S$ is
bijective.\\

(iii) We have already seen that $S$ is an
anti-algebra
automorphism of $H$ and ${\rm dim}_k\int _l^H=1$. The rest of the
proof is
identical to the proof for classical Hopf
algebras.
\end{proof}

\begin{remark}\rm 
We cannot deduce the isomorphism
$\theta $ in Theorem 2.2 from a Structure 
Theorem for dual quasi-Hopf
bicomodules. If $A$ is a dual quasi-Hopf algebra 
(for the complete
definition see the last Section) then the category of $A$-bicomodules
$^A{\cal M}^A$ is monoidal and $A$ is in a canonical way an algebra in
$^A{\cal M}^A$. 
Thus, it makes sense to define a right dual quasi-Hopf
$A$-bicomodule $M$ as being a 
right $A$-module in $^A{\cal M}^A$. Denote
by $^A{\cal M}^A_A$ the category whose 
objects are right dual quasi-Hopf
$A$-bicomodules and morphisms $A$-bicomodule maps which are also right
$A$-linear (for more details see \cite{sh}). This definition is dual to the
one given by Hausser and Nill \cite{hn3} for quasi-Hopf bimodules. So,
using their 
Structure Theorem, by duality, we can prove a Structure
Theorem for dual quasi-Hopf 
$A$-bicomodules and then we can apply it in
the  particular case $M=H$ (here $H$ 
is a finite dimensional quasi-Hopf
algebra and $A=H^*$, the linear dual 
space of $H$). But, on a hand, in
order to obtain $H$ as an object in $^{H^*}{\cal M}^{H^*}_{H^*}$ we need
the antipode $S$ of $H$ to be bijective and, on the other hand, using the
definition of Hausser and Nill for coinvariants, in the dual case we don't
obtain the space of integrals in $H$. In consequence, if the isomorphism
$\theta $ is derived  from a Structure Theorem for dual quasi-Hopf
bicomodules then it is not the dual case of Hausser and Nill result.\\
Now, let $M\in ^A{\cal M}^A_A$. First, since $A$ is a coassociative
coalgebra  we can define the set of right coinvariants $M^{co(A)}$ of $M$
as  being $M^{co(A)}:=\{m\in M\mid \r _r(m)=m\ot 1\}$, 
where $\r _r:M\ra M\ot A$, $\r _r(m)=\sum m_{(0)}\ot m_{(1)}$, $m\in M$, is the 
right coaction of $A$ on $M$. Since $\r _r(m\cd a)=\r _r(m)\Delta (a)$ for all
$m\in M$ and $a\in A$ (we denote by $M\ot A\ni m\ot A\ra m\cd a\in M$,
$m\in M$, $a\in A$, the right action of $A$ on $M$), as in \cite[Theorem
2.11]{bn} we can show that $\ov {P}: M\ra M^{co(A)}$ given by $\ov
{P}(m):=\sum m_{(0)}\cd \b (m_{(1)})S(m_{(2)})$, $m\in M$, is well defined
and a surjection. Thus, the map $\ov {\theta }^{-1}:M\ra M^{co(A)}\ot A$,
$\ov {\theta }^{-1}(m):=\sum \ov {P}(m_{(0)})\ot m_{(1)}$, $m\in M$, is
well defined. It is not hard to see that $M^{co(A)}\ot A\in ^A{\cal M}^A_A$
via the structures $(m\ot a)\cd b:=m\ot ab$, and $M^{co(A)}\ot A\ni m\ot a
\ra \sum a_1\ot m\ot a_2\ot a_3\in A\ot M^{co(A)}\ot A\ot A$, $m\in
M^{co(A)}$, $a, b\in A$. In this way $\ov {\theta }^{-1}$ becomes a right
$A$-colinear map but not a morphism in $^A{\cal M}^A_A$; moreover, we don't
know if, in general, $\ov {\theta }^{-1}$ is bijective. Under these
circumstances, if we take $H$ a finite dimensional quasi-Hopf algebra then
$H\in ^{H^*}{\cal M}^{H^*}$ via 
$$
H\ni h\ra \sum \limits _{i,
j=1}^ne^i\ot e_jhe_i\ot e^j\in H^*\ot H\ot H^*,
$$
where $\{e_i\}_{i=\ov
{1, n}}$ is a basis of $H$ with dual basis $\{e^i\}_{i=\ov {1, n}}$ 
of $H^*$. Moreover, if we define 
$$
h\cd h^*:=\sum h^*(S(X^2h_2)\a X^3)X^1h_1
\mbox{${\;\;\;}$$\forall $ $h\in H$, $h^*\in H^*$}
$$
then with
this structures $H$ is not an object in $^{H^*}{\cal M}^{H^*}_{H^*}$ but
$\r _r(h\cd h^*)=\r _r(h)\tilde {\Delta }(h^*)$, $H^{co(H^*)}=\int _l^H$
and the projection $\ov {P}: H\ra \int _l^H$ is just the projection $P$
defined in (\ref{p}). Therefore, in this case $\ov {\theta }^{-1}$ coincide
with $\theta ^{-1}$ defined in (\ref{thi}), so it is bijective. We notice
that, in the  proof of the fact that $\theta ^{-1}$ defined by (\ref{thi})
is the inverse of $\theta $, 
a key role is played by the relation
(\ref{f2}) in Lemma 2.1. Since the equality (\ref{f2}) involve all the
structures of $H$ as a quasi-Hopf algebra, it cannot be generalized for the
coinvariants $M^{co(A)}$ of a dual quasi-Hopf $A$-bicomodule $M$. In
conclusion, in order to obtain $\theta $ from a Structure Theorem for dual
quasi-Hopf bicomodules we need the above context but it don't provide a
"suggestion" for the general case. Also, the same kind of problems occur
when we work with relative Hopf modules \cite{bn}, instead of dual
quasi-Hopf bicomodules.\\          
\end{remark}

Let $H$ be a quasi-Hopf
algebra and $t$ a left integral
in $H$. Using the fact that $H$ is an
associative algebra, we find that
$th$ is also a left integral in $H$, for
all $h\in H$, hence the
space of left (right) integral in $H$ is a
two-sided ideal.
Moreover, if $H$ is finite dimensional, then it follows
from the uniqueness of
the integral in $H$, that there exists $\mu \in H^*$
such that
\begin{equation}\label{gdi} th=\mu
(h)t,
\mbox{${\;\;\;}$
$\forall $ $t\in \int _l^H$ and $h\in
H$.}
\end{equation}
More precisely, $\mu \in {\rm Alg}(H, k)$. It was noted
in \cite{hn3}
that ${\rm Alg}(H, k)$ is a group with multiplication given
by
$\nu \xi =(\nu \ot \xi )\circ \Delta $, unit $\varepsilon$,
and inverse
$\mu ^{-1}=\mu \circ S$. In \cite{hn3}, $\mu$ is
called the modulus of $H$,
but, following the classical terminology
for Hopf algebras, we will call
$\mu $ the distinguished grouplike element.
Observe that $\mu =\va $ if and
only if $H$ is unimodular. As
in the case of a Hopf algebra, it follows
from the bijectivity of
the antipode that
\begin{equation}\label{gdim}
hr=\mu ^{-1}(h)r= \mu (S(h))r,
\mbox{${\;\;\;}$$\forall $ $r\in \int _r^H$
and $h\in H$.}
\end{equation}
For infinite dimensional Hopf algebras it is
well known
that the space of left (right) integrals in $H$ is
zero \cite[pg. 107]{sw}. In order to prove a similar result for
quasi-Hopf algebras we first need a Lemma.

\begin{lemma}\lelabel{2.3}
Let $H$ be a
quasi-Hopf algebra in the sense of Drinfeld and
define
\begin{equation}
\un
{\Delta }: H\ra H\ot H,
\mbox{${\;\;}$}
\un {\Delta }(h)=\sum h_{\un
{1}}\ot h_{\un {2}}:=\sum
q^1h_1p^1\ot q^2h_2p^2,
\mbox{${\;\;\;}$$\forall
h\in H$,}
\end{equation}
where $p_R=\sum p^1\ot p^2$ and $q_R=\sum q^1\ot
q^2$ are defined
by (\ref{qr}). If $J$ is a non-zero two-sided ideal of $H$
such
that $\un {\Delta }(J)\subseteq J\ot H$, then
$J=H$.
\end{lemma}

\begin{proof}
From (\ref{pqr}), we easily deduce
that
$$
\sum (1\ot S^{-1}(p^2))\un {\Delta }(p^1hq^1)(1\ot S(q^2))
=\Delta
(h),
\mbox{${\;\;\;}$$\forall h\in H$.}
$$
This implies that $\Delta
(J)\subseteq
J\ot H$,
since $J$ is a two-sided ideal of $H$ and $\un
{\Delta }(J)\subseteq
J\ot H$. Now, if $\va (J)=0$, then for any $h\in H$
we have
$h=\sum \va (h_1)h_2\in \va (J)H=0$, so $J=0$, a
contradiction.
Thus $\va (J)\not=0$, and there exists $a\in J$ with $\va
(a)=1$.
Using (\ref{q5}), we obtain that
$\b =\va (a)\b =\sum a_1\b
S(a_2)\in JH\subseteq J$, so $\b \in J$.
Using (\ref{q6}) and the fact that
$J$ is a two-sided ideal of $H$,
we find that $1=\sum X^1\b S(X^2)\a X^3\in
J$, and $J=H$.
\end{proof}

For $h\in H$ and $h^*\in H^*$, we
define
$h^*\act h=\sum h^*(h_{\un {2}})h_{\un {1}}$.
For a two-sided ideal
$I$ of $H$, we let $H^*\act I$
be the subspace of $I$ generated by all the
elements of the form
$h^*\act a$, with $h^*\in H^*$ and $a\in
I$.

\begin{theorem}\thlabel{2.4}
Let $H$ be a quasi-Hopf algebra in the
sense of Drinfeld and $I$ a
non-zero two-sided ideal of $H$.
Then
$$
J=H^*\ra I=H.
$$
As a consequence, we obtain\\
${\;\;}$(i)
If $H$
is a quasi-Hopf algebra with an
antipode $S$, then $H$ is finite
dimensional if and only if $S$ is
bijective and $\int
_l^H\not=0$.\\
${\;\;}$(ii)
If a quasi-Hopf algebra in the sense of
Drinfeld is semisimple as
an algebra, then it is
finite
dimensional.
\end{theorem}

\begin{proof}
The statement follows from
\leref{2.3} if we can show that $J$ is a non-zero
two-sided ideal of $H$
such that $\un {\Delta }(J)\subseteq J\ot H$.\\
Obviously $\va \act h=h$,
and therefore
$I\subseteq J$. For all $h\in H$, $h^*\in H^*$ and $a\in I$
we have
\begin{eqnarray*}
(h^*\ra a)h
&=&\sum
h^*(q^2a_2p^2)q^1a_1p^1h\\
{\rm (\ref{qr1})}~~~~ &=&\sum
h^*(q^2(ah_1)_2p^2S(h_2))q^1(ah_1)_1p^1\\
&=&\sum (S(h_2)\rh h^*)\ra (ah_1)
\end{eqnarray*}
and $J$ is a right ideal. $J$ is also a left ideal,
since
\begin{eqnarray*}
h(h^*\ra a)%
&=&\sum
h^*(q^2a_2p^2)hq^1a_1p^1\\
{\rm (\ref{qr1})}~~~~&=&\sum
h^*(S^{-1}(h_2)q^2(h_1a)_2p^2)q^1(h_1a)_1p^1\\
&=&\sum (h^*\lh
S^{-1}(h_2))\ra (h_1a).
\end{eqnarray*}
Write $f=\sum f^1\ot f^2$. Using
(\ref{qr2}), (\ref{q1}) and (\ref{pr}) we can show 
that\\[2mm]

${\;\;\;}$
$h^*\act (g^*\act h)$
$$
=\sum [(g^1S(x^3)\rh h^*\lh
S^{-1}(f^2X^3))(g^2S(x^2)\rh g^*\lh
S^{-1}(f^1X^2))]\act (X^1ax^1)
$$
for all $h^*, g^*\in H^*$ and $h\in H$. $I$ is a two-sided ideal 
of $H$, so the
above equality shows
that $H^*\ra J\subseteq J$. To prove that $\un
{\Delta
}(J)\subseteq J\ot H$, we use the same arguments as in
\cite[p. 12]{m}.
Take $a\in J$, and write $\un {\Delta }(a) =\sum
\limits _{i=1}^n a_i\ot
a'_i$, where $a_1, \cdots , a_m\in J$
and $a_{m + 1}, \cdots , a_n$ are
linearly independent modulo $J$.
For any $h^*\in H^*$, $h^*\act a=\sum
\limits
_{i=1}^nh^*(a'_i)a_i\in J$. The linear independence of $a_{m +
1},
\cdots , a_n$ modulo $J$ implies that $h^*(a'_i)=0$, and
therefore
$a'_i=0$ ($h^*$ is arbitrary),
for all $i>m$. We find that $\un {\Delta
}(a)\in
J\ot H$, as needed.\\

(i) One implication follows from
\thref{2.2}. Conversely, assume that
$S$ is bijective and $I=\int
_l^H\neq
0$. Then $I$ is a non-zero two-sided ideal of $H$ and the first
part of our
Theorem tells us that $H^*\act I=H$. Thus there exist
$\{h^*_i\}_{i=\ov {1,
n}}\subseteq H^*$ and $\{t_i\}_{i=\ov {1,
n}} \subseteq \int _l^H$ such
that $1=\sum \limits
_{i=1}^nh^*_i\act t_i$. For any $i=\ov {1, n}$ we have
$\un
{\Delta }(t_i)=\sum \limits _{j=1}^{n_i}a^i_j\ot b^i_j$, for
some
$\{a^i_j\}_{j=\ov {1, n_i}}\subseteq H$ and $\{b^i_j\}_{j=\ov
{1,
n_i}}\subseteq H$. Therefore, for any $h^*\in H^*$ and $i=\ov {1,
n}$
we have $h^*\ra t_i=\sum \limits _{j=1}^{n_i}
h^*(b^i_j)a^i_j$. For
all
$h\in H$ we obtain that
\begin{eqnarray*}
h
&=&\sum \limits
_{i=1}^nh(h^*_i\ra t_i)\\
&=&\sum \limits _{i=1}^n(h^*\lh S^{-1}(h_2))\ra
h_1t_i\\
&=&\sum \limits _{i=1}^n(h^*\lh S^{-1}(h))\ra
t_i
\mbox{${\;\;\;}$since $t_i\in \int _l^H, \forall i=\ov {1,
n}$}\\
&=&\sum \limits _{i=1}^n\sum \limits
_{j=1}^{n_i}
h^*(S^{-1}(h)b^i_j)a^i_j.
\end{eqnarray*}
This shows that $H$
is finite dimensional, since it is a subspace of the
span of $\{a^i_j\mid
i=\ov {1, n}, j=\ov {1, n_i}\}$.\\

(ii) Let $H$ be a semisimple
quasi-Hopf
algebra with bijective antipode. Then ${\rm Ker}(\va )$ is a
two-sided
ideal of $H$. Since $H$ is a semisimple left $H$-module,
there exists a
left ideal $I$ of $H$ such that $H=I\oplus {\rm Ker}(\va
)$. ${\rm Ker}(\va
)$ has codimension $1$ in $H$, hence
$I$ has dimension $1$. Write $1=t +
h$, with $t\in I$, $h\in {\rm Ker}(\va
)$. $t\not=0$, because $1\not\in
{\rm Ker}(\va )$.
It follows that $I=kt$, since ${\rm dim}_k(I)=1$.
For all
$h'\in H$, we have $h't\in I$, and also
$h't=\va (h')t +
(h' - \va (h')1)t$
with $\va (h')t\in I$ and $(h' - \va
(h')1)t\in {\rm Ker}(\va )$. Since we
have a direct sum, it follows that
$h't=\va (h')t$, and $t$ is a non-zero
left integral in $H$. 
From (i), it then follows that $H$ is finite
dimensional.
\end{proof}

\begin{remarks}\relabel{2.5}\rm
(i)
Let $H$ be a
quasi-Hopf algebra in the sense of Drinfeld. Then
$S(\int _l^H)=\int _r^H$
and ${\rm dim}_k\int _l^H=
{\rm dim}_k\int _r^H$. Therefore, if $H$ is
infinite dimensional then
$\int _l^H=\int _r^H=0$.\\

(ii)
Let $H$ be a
finite dimensional quasi-Hopf algebra
and $t$ a non-zero left integral in
$H$.
\thref{2.4} implies that $H=H^*\act t$ so the map
\begin{equation}
\ov
{\theta }:\ H^*\ra H,
\mbox{${\;\;}$}
\ov {\theta }(h^*)=h^*\ra t=\sum
h^*(q^2t_2p^2)q^1t_1p^1,
\mbox{${\;}$$\forall h^*\in
H^*$}
\end{equation}
is bijective. Moreover, $\ov {\theta }$ is a left
$H$-linear,
where the left $H$-action on $H^*$ is given by the
formula
$$h\cd h^*=h^*\lh S^{-1}(h)$$
for $h\in H$ and $h^*\in H^*$.
If $H$
is a classical Hopf algebra, then $\ov {\theta }$
is also left $H^*$-linear
where the $H^*$-actions on $H^*$ and $H$
are given by convolution
and
$h^*\rh h=\sum h^*(h_2)h_1$ respectively. This means that $H$ is a
left
cyclic $H^*$-module generated by a left non-zero
integral.
\end{remarks}

\section{Cointegrals on quasi-Hopf
algebras}\selabel{3}
\setcounter{equation}{0}
In the first part of this
Section we study the cointegrals on
a finite dimensional quasi-Hopf algebra
(so the
antipode of $H$ is bijective).

\begin{definition}\delabel{3.1}
\cite{hn3}
Let $H$ be a quasi-Hopf algebra, $M$ an $H$-bimodule and
$\rho
:\ M\ra M\ot H$ an $H$-bimodule map. Then $(M, \rho )$
is called a right
quasi-Hopf $H$-bimodule if thefollowing
relations
hold:
\begin{equation}\label{hb1}
(id_M\ot \va )\circ \rho
=id_M,
\end{equation}
\begin{equation}\label{hb2}
\Phi\cd (\rho \ot
id_M)(\rho (m))=(id_M\ot \Delta )(\rho (m))\cd \Phi
,
\mbox{${\;\;\;}$$\forall m\in M$.}
\end{equation}
\end{definition}

A morphism between two right quasi-Hopf $H$-bimodules is an
$H$-bimodule map
$f:\ M\ra M'$ satisfying $\rho'\circ f=(f\ot id)\circ \rho $.
$_H{\cal
M}_H^H$ is the category of right quasi-Hopf $H$-bimodules
and morphisms of
right quasi-Hopf $H$-bimodules.\\
We will use the Sweedler type notation
$$
\rho (m)=\sum m_{(0)}\ot m_{(1)},
\mbox{${\;\;}$} (\rho \ot id_M)(\rho
(m))=m_{(0, 0)}\ot m_{(0,
1)}\ot m_{(1)}, \mbox{${\;\;}$etc.} $$
Let $H$ be
a
quasi-Hopf algebra and $M$ a right quasi-Hopf $H$-bimodule. We define
$E:
M\ra M$,
by
\begin{equation}\label{eg}
E(m)=\sum X^1\cd m_{(0)}\cd \b
S(X^2m_{(1)})\a X^3,
\end{equation}
for all $m\in M$.
$$M^{{\rm
co}H}=\{n\in M\mid E(n)=n\}$$
is called the space of
coinvariants of $M$.
We also have (cf. \cite[Corollary
3.9]{hn3})
\begin{equation}\label{coinv}
M^{{\rm co}H}=\{n\in M\mid \rho
(n)=\sum E(x^1\cd
n)\cd x^2\ot x^3\}.
\end{equation}
We have the following
structure Theorem
for right quasi-Hopf $H$-bimodules (see \cite[Theorem
3.8]{hn3}).
The map
\begin{equation}
\nu :\ M^{{\rm co}H}\ot H\ra
M,
\mbox{${\;}$}
\nu (n\ot h)=n\cd h,
\mbox{${\;}$$\forall n\in M^{{\rm
co}H},~ h\in H$}
\end{equation}
is an isomorphism of right quasi-Hopf
$H$-bimodules. Here
$M^{{\rm co}H}\ot H$ is a right  quasi-Hopf
$H$-bimodule via the structures
$a\cd (n\ot h)\cd b=\sum E(a_1\cd n)\ot
a_2hb$ and
$\rho (n\ot h)=\sum E(x^1\cd n)\ot x^2h_1\ot x^3h_2$, $\forall
n\in N,
a, h, b\in H$. The inverse of $\nu $ is given
by
\begin{equation}
\nu ^{-1}(m)=\sum E(m_{(0)})\ot
m_{(1)},
\mbox{${\;\;\;}$$\forall m\in M$.}
\end{equation}
Now, let $H$ be
a finite dimensional quasi-Hopf algebra; recall that
the antipode is then
automatically bijective. Let
$\{e_i\}_{i=\ov {1, n}}$ be a basis in $H$
with dual basis
$\{e^i\}_{i=\ov {1, n}}$ in $H^*$ and
consider
\begin{equation}\label{uv}
U=\sum g^1S(q^2)\ot
g^2S(q^1),
\mbox{${\;\;\;}$}
V=\sum S^{-1}(f^2p^2)\ot
S^{-1}(f^1p^1),
\end{equation}
where $f=\sum f^1\ot f^2$ is the element
defined by (\ref{f}),
$f^{-1}=\sum g^1\ot g^2$, and $q_R=\sum q^1\ot q^2$
and
$p_R=\sum p^1\ot p^2$ are defined as in (\ref{qr}).
Following
\cite{hn3}, $H^*$ is right quasi-Hopf $H$-bimodule.
The structure is the
following
\begin{equation}\label{str1}
h\cd h^*\cd h'=S(h')\rh h^*\lh
S^{-1}(h),
\mbox{${\;\;\;}$$\forall h, h'\in H, h^*\in
H^*$,}
\end{equation}
\begin{equation}\label{str2}
\rho (h^*)=\sum \limits
_{i=1}^ne^i * h^*\ot e_i,
\mbox{${\;\;\;}$$\forall h^*\in
H^*$,}
\end{equation}
where the (non-associative) multiplication $* :\
H^*\ot H^*\ra H^*$
is given by
\begin{equation}
\langle h^* * g^*, h\rangle
=\sum h^*(V^1h_1U^1)g^*(V^2h_2U^2),
\mbox{${\;\;\;}$$\forall h^*, g^*\in
H^*, h\in H$.}
\end{equation}
The coinvariants $\l \in H^{*{\rm co}H}$ are
called left cointegrals on $H$
and the space of left cointegrals is denoted
by ${\cal L}$. Thus
${\rm dim}_k{\cal L}=1$ and the projection $E:\ H^*\ra
{\cal L}$ is given by
\begin{equation}\label{e}
\langle E(h^*), h\rangle
=\sum \limits _{i=1}^n\langle e^i\ot h^*, \ov
{\Delta
}(S^{-1}(q^1)hS^2(q^2e_i)S(\b ))\rangle ,%
\end{equation}
for all
$h^*\in
H^*$ and $h\in H$.
Here $\ov {\Delta }(h)=\sum V^1h_1U^1\ot V^2h_2U^2$, for
any
$h\in H$. The transpose $E^T:\ H\ra H$ is given
by
\begin{equation}\label{et}
E^T(h)=\sum \limits _{i=1}^n\langle e^i\ot
id, \ov {\Delta}
(S^{-1}(q^1)h S^2(q^2e_i)S(\b ))\rangle
,%
\mbox{${\;\;}$$\forall h\in H$}
\end{equation}
and provides a projection
onto the space of right integrals
in $H$. Moreover, the dual pairing
${\cal
L}\ot \int _r^H\to k$, $\lambda\ot r\mapsto \langle \l , r\rangle $
is
nondegenerate, \cite[Lemma 4.4]{hn3}.\\
When $H$ is an ordinary Hopf
algebra the
left cointegrals on $H$ are precisely the left integrals in
$H^*$
(i.e., an element $\l \in H^*$ such that $\sum \l (h_2)h_1=\l
(h)1,
\forall h\in H$). For quasi-Hopf algebras $H$ we will give
some
characterizations for left cointegrals.
First, we need another formula for
the projection (\ref{e}),
giving the connection between the projection $P$

onto the space of left integrals defined by (\ref{p}) and
the projection
$E^T$, defined by (\ref{et}). This provides an
answer to a question raised
in \cite{hn3}.

\begin{lemma}\lelabel{3.2}
Let $H$ be a finite dimensional
quasi-Hopf algebra and $E$ the map defined
by (\ref{e}). Then, for all $h^*
\in H^*$ and $h\in H$, we have:
\begin{equation}\label{se}
\langle E(h^*),
h\rangle =\langle h^* , S^{-1}(P(S(h)))\rangle .
\end{equation}
In
particular, if $E^T$ is the transpose map (\ref{et}) and
$P: H\ra \int
_l^H$ is the projection (\ref{p}) then
$E^T(h)=S^{-1}(P(S(h)))$, for all
$h\in H$.
\end{lemma}

\begin{proof}
Take $f=\sum f^1\ot f^2=\sum F^1\ot
F^2$ and $f^{-1}=\sum g^1\ot g^2
=\sum G^1\ot G^2$
as in (\ref{f}) and
(\ref{g}). It is easy to see that
\begin{equation}\label{fgab}
\sum
g^1S(g^2\a )=\b ,
\mbox{${\;\;\;}$}
\sum S(\b f^1)f^2=\a
,
\mbox{${\;\;\;}$}
\sum f^1\b S(f^2)=S(\a )
\end{equation}
and we compute,
for $h^*\in H^*$ and $h\in H$:
\begin{eqnarray*}
\hspace*{-2cm} \langle
E(h^*), h\rangle &=&  
\sum \limits _{i=1}^n
\langle e^i\ot h^*, \ov
{\Delta }(S^{-1}(q^1)hS^2(q^2e_i)S(\b ))\rangle \\
{\rm (\ref{uv},
\ref{ca})}~~~~&=&\sum \limits _{j=1}^n
\langle e^j,
hS^2(q^2S^{-1}(f^2q^1_2p^2)(e_j)_1U^1)S(\b )\rangle
\langle h^*,
S^{-1}(f^1q^1_1p^1)(e_j)_2U^2\rangle \\
{\rm
(\ref{pqr},\ref{uv},\ref{qr})}~~~~&=&\sum \limits _{j=1}^n
\langle e^j,
hS^2(S^{-1}(f^2)(e_j)_1g^1S(X^2)\a X^3S^{-1}(\b ))\rangle \\
&&~~~~\langle
h^*, S^{-1}(f^1)(e_j)_2g^2S(X^1)\rangle \\
{\rm
(\ref{g2},\ref{pf})}~~~~&=&\sum \limits _{j=1}^n
\langle e^j,
hS^2(S^{-1}(f^2)(e_j)_1F^2_1X^2g^1_2G^2\a
S^{-1}(\b F^1X^1g^1_1G^1))\rangle
\\
&&~~~~\langle h^*, S^{-1}(f^1)(e_j)_2F^2_2X^3g^2\rangle \\
{\rm
(\ref{fgab}, \ref{q5})}~~~~&=&\sum \limits _{i=1}^n
\langle e^i,
hS^2(S^{-1}(f^2)(e_i)_1X^2)S(\b F^1X^1\b )F^2\rangle
\langle h^*,
S^{-1}(f^1)(e_i)_2X^3\rangle \\
{\rm (\ref{fgab},\ref{ca})}~~~~&=&\sum
\limits _{i=1}^n
\langle e^i,
hS^2(S^{-1}(S(e_i)_2)S^{-1}(f^2)X^2S^{-1}(S^{-1}(\a )X^1\b
))\rangle
\\
&&~~~~\langle h^*, S^{-1}(S(e_i)_1)S^{-1}(f^1)X^3\rangle \\
{\rm
(\ref{g2},\ref{pf})}~~~~&=&\sum \limits _{i=1}^n
\langle e^i,
hS(S(e_i)_2)S(\b F^1f^2_1X^2g^1_2)F^2f^2_2X^3g^2\a \rangle \\
&&~~~~\langle
h^*, S^{-1}(S(e_i)_1)S^{-1}(f^1X^1g^1_1)\rangle \\
{\rm (\ref{fgab},
\ref{q5})}~~~~&=&\sum \limits _{j=1}^n
\langle e^j,
g^1S(hS((e_j)_2)S(X^2)\a X^3g^2\a )\rangle \langle
h^*,
S^{-1}(X^1(e_j)_1)\rangle \\
{\rm (\ref{fgab}, \ref{qr},
\ref{p})}~~~~&=&\langle h^*, S^{-1}(P(S(h)))\rangle .
\end{eqnarray*}
The
last assertion follows easily by (\ref{se}) and
(\ref{et}).
\end{proof}

\begin{lemma}\lelabel{3.3}
Let $H$ be a finite
dimensional quasi-Hopf algebra and
$E:\ H^*\ra  {\cal L}$ the projection
(\ref{e}). Then
\begin{equation} \label{f3}
E(h^*\lh S^{-1}(h))=\mu
(h)E(h^*),
\mbox{${\;\;\;}$$\forall h^*\in H^*$ and $h\in
H$,}
\end{equation}
\begin{equation}  \label{f4}
\l (S^{-1}(h)h')=\sum \mu
(h_1)\l (h'S(h_2)),
\mbox{${\;\;\;}$$\forall h, h'\in
H$,}
\end{equation}
where $\mu $ is the distinguished grouplike element of
$H^*$.
\end{lemma}

\begin{proof}
The relations (\ref{f3}) follow easily
from (\ref{se})
and the fact that $P(h)\in \int _l^H$, $\forall h\in H$. To
prove (\ref{f4}),
let $\{e_i\}_{i=\ov {1, n}}$ be a basis in $H$ with dual
basis
$\{e^i\}_{i=\ov {1, n}}$ in $H^*$ and
$q_R$ the element defined by
(\ref{qr}).
For all $h, h'\in H$ and $\l \in {\cal L}$ we
have:
\begin{eqnarray*}
&&\hspace*{-2cm}\l (S^{-1}(h)h')=
\langle E(\l ),
S^{-1}(h)h'\rangle \\
{\rm (\ref{se},\ref{p})}~~~~&=&\sum \limits
_{i=1}^n\langle e^i, \b
S^2(q^2(e_i)_2)S(h')h\rangle
\langle \l ,
S^{-1}(q^1(e_i)_1)\rangle \\
&=&\sum \limits _{j=1}^n\langle e^j, \b
S^2(q^2(e_j)_2h_2)S(h')\rangle
\langle \l , S^{-1}(q^1(e_j)_1h_1)\rangle
\\
{\rm (\ref{p},\ref{se})}~~~~&=&\sum \langle E(\l \lh
S^{-1}(h_1)),
h'S(h_2)\rangle \\
{\rm (\ref{f3})}~~~~&=&\sum \mu (h_1)\l
(h'S(h_2)).
\end{eqnarray*}
\end{proof}

\begin{proposition}\prlabel{3.4}
Let $H$ be a finite dimensional quasi-Hopf algebra and $\mu $ th
e
distinguished grouplike element of $H^*$. For $\lambda\in H^*$,
the
following
statements are equivalent:
\begin{enumerate}
\item[(a)] $\l$ is a
left cointegral on $H$;
\item[(b)] for all $h\in H$, we
have
\begin{eqnarray}
&&\hspace*{-2cm}\sum
\l
(S^{-1}(f^1)h_2S^{-1}(q^1g^1))S^{-1}(f^2)h_1S^{-1}(q^2g^2)
\nonumber\\
&=&\sum \mu (q^1_1x^1)\langle \l ,
hS^{-1}(f^1)g^2S(q^1_2x^2)\rangle
q^2x^3S^{-1}(S^{-1}(f^2)g^1);\label{coin}
\end{eqnarray}
\item[(c)] for all
$h\in H$, we have
$$\sum \l (\smi (f^1)h_2U^2)\smi (f^2)h_1U^1=
\sum \mu
(q^1_1x^1)\langle \l , hS(q^1_2x^2)\rangle
q^2x^3.\label{co2}
$$
\end{enumerate}
Here $f=\sum f^1\ot f^2$, $q_R=\sum
q^1\ot q^2$ and $U=\sum U^1\ot U^2$ are
defined
respectively  by (\ref{f}),
(\ref{qr}) and (\ref{uv}), and
$f^{-1}=\sum g^1\ot
g^2$.
\end{proposition}

\begin{proof}
$\un{(a)\Longrightarrow
(b)}$.
Suppose that $\l $ is a left cointegral. As before, we write
$f=\sum
f^1\ot f^2=\sum F^1\ot F^2=\sum {\bf F}^1\ot {\bf F}^2,
f^{-1}=\sum g^1\ot
g^2=\sum G^1\ot G^2=\sum {\bf G}^1\ot {\bf G}^2$,
$q_R=\sum q^1\ot q^2$ and
$p_R=\sum p^1\ot p^2$. Using (\ref{coinv}),
(\ref{str1}),
(\ref{str2}) and
(\ref{f3}), we find that
\begin{equation} \label{cointe}
\sum \l
(V^2h_2U^2)V^1h_1U^1=\sum \mu (x^1)\l (hS(x^2))x^3,
\end{equation}
for all
$h\in H$, and we compute that
\begin{eqnarray*}
&&\hspace*{-2cm}\sum \mu
(q^1_1x^1)\langle \l,
hS^{-1}(f^1)g^2S(q^1_2x^2)\rangle
q^2x^3S^{-1}(S^{-1}(f^2)g^1)\\
{\rm(\ref{f3}, \ref{cointe})}~~~~&=&\sum
\langle \l
,
V^2[S^{-1}(q^1)hS^{-1}(f^1)g^2]_2U^2\rangle
q^2V^1[S^{-1}(q^1)hS^{-1}(f^1)g^2]
_1U^1\\%
&&~~~~S^{-1}(S^{-1}(f^2)g^1)\\
{\rm(\ref{uv},\ref{ca})}~~~~&=&\sum
\langle \l , S^{-1}({\bf F}^1q^1_1p^1)h_2
S^{-1}(F^1f^1_1{\bf
G}^1)g^2_2U^2\rangle q^2S^{-1}({\bf
F}^2q^1_2p^2)h_1
S^{-1}(F^2f^1_2{\bf
G}^2)\\
&&~~~~g^2_1U^1S^{-1}(S^{-1}(f^2)g^1
)\\
{\rm(\ref{pqr},\ref{ca})}~~~~&=&\sum \langle \l
,
S^{-1}(S(h)_1)S^{-1}(F^1f^1_1)g^2_2G^2S(X^1)\rangle
S^{-1}(S(h)_2)S^{-1}(F^2f^1_2)\\
&&~~~~g^2_1G^1S(X^2)\a S^{-1}(S^{-1}(f^2)g^1S(X^3))\\
{\rm(\ref{g2},\ref{pf})}~~~~&=&\sum \langle \l , 
S^{-1}(S(h)_1)S^{-1}(S(X^3)F^1f^1_1)g^2\rangle
S^{-1}(S(h)_2)S^{-1}(S(X^2)F^2f^1_2)\\
&&~~~~g^1_2G^2\a
S^{-1}(S^{-1}(S(X^1)f^2)g^1_1G^1)\\
{\rm(\ref{fgab},\ref{q5})}~~~~&=&\sum
\langle \l
,
S^{-1}(S(h)_1)S^{-1}(S(X^3)F^1f^1_1)\rangle
S^{-1}(S(h)_2)\\
&&~~~~S^{-1}(S^{-
1}(S(X^1)f^2)\b S(X^2)F^2f^1_2)\\
{\rm(\ref{g2},\ref{pf})}~~~~&=&\sum
\langle \l , 
S^{-1}(S(h)_1)S^{-1}(f^1X^1)\rangle
S^{-1}(S(h)_2)S^{-1}(S^{-1}(F^2f^2_2X^3)\b
F^1f^2_1X^2)\\
{\rm(\ref{fgab},\ref{q5})}~~~~&=&\sum \langle \l ,
S^{-1}(X^1S(h)_1)\rangle
S^{-1}(S^{-1}(\a X^3)X^2S(h)_2)\\
{\rm(\ref{qr},\ref{ca})}~~~~&=&\sum \langle \l
,
S^{-1}(f^1)h_2S^{-1}(q^1g^1)\rangle
S^{-1}(f^2)h_1S^{-1}(q^2g^2).
\end{eqnarray*}
$\un{(b)\Longrightarrow (a)}$.
Assume that $\l \in H^*$ satisfies
(\ref{coin}).
It follows from (\ref{ca})  that
\begin{equation}
\sum \l
(S^{-1}(q^1h_1))q^2h_2=
\sum \mu (q^1_1x^1)
\langle \l ,
S^{-1}(f^1h)g^2S(q^1_2x^2)\rangle
S^{-1}(f^2)g^1S(q^2x^3)\label{scoin},
\end{equation}
for all $h\in H$ and
\begin{eqnarray*}
&&\hspace*{-2cm}\langle \l ,
S^{-1}(P(S(h)))\rangle \\
{\rm(\ref{p}, \ref{scoin})}~~~~&=&\sum \mu
(q^1_1x^1)\langle \l ,
S^{-1}(f^1\b
S(f^2)S^2(g^1)S^3(q^2x^3)S(h))g^2S(q^1_2x^2)\rangle
\\
{\rm(\ref{fgab})}~~~~&=&\sum \mu (q^1_1x^1)\langle \l, 
hS^2(q^2x^3)S(g^1)\a g^2S(q^1_2x^2)\rangle \\
{\rm(\ref{fgab}, \ref{qr},
\ref{qr1})}~~~~&=&\l (h),
\end{eqnarray*}
and it follows from \leref{3.2}
that $\l \in {\cal L}$.\\
$\un{(c)\Longrightarrow (a)}$.
Repeating the
computations of the first part of the proof of
\leref{3.2}, we find that
the projection $E$ (cf. (\ref{e})) is given by
\begin{equation}
\label{the}
\langle E(h^*), h\rangle =\sum \limits _{i=1}^n \langle
e^i, hS(f^2)S^2((e_i)_1U^1)S(\b )\rangle
\langle h^*,
S^{-1}(f^1)(e_i)_2U^2\rangle ,
\end{equation}
for all $h^*\in H$, $h\in H$.
Using (\ref{the}), we can compute that
$E(\lambda)=\lambda$, so $\l \in
{\cal L}$.\\
$\un{(a)\Longrightarrow (c)}$. Assume that $\l \in {\cal
L}$.
We calculate
\begin{eqnarray*}
&&\hspace*{-2cm}\sum \mu (q^1_1x^1)\langle
\l , hS(q^1_2x^2)\rangle q^2x^3\\
{\rm(\ref{f4}, \ref{cointe})}~~~~&=&\sum
\langle \l ,
V^2[S^{-1}(q^1)h]_2U^2\rangle
q^2V^1[S^{-1}(q^1)h]_1U^1\\
{\rm(\ref{uv},
\ref{ca}, \ref{pqr})}~~~~&=&
\sum \langle \l , S^{-1}(f^1)h_2U^2\rangle 
S^{-1}(f^2)h_1U^1
\end{eqnarray*}
and the proof is
complete.
\end{proof}

\begin{remark}\rm
Formula (\ref{scoin}) is
equivalent to (\ref{coin}), and can be
viewed as another characterization
of left cointegrals. In the
case of a Hopf algebra, (\ref{scoin}) takes the
form
$$\sum \lambda(S^{-1}(h_1))h_2=\lambda(S^{-1}(h))1$$
which is the
well-known statement that $\lambda$ is a left
integral if and only if
$\lambda\circ S^{-1}$ is a right
integral.
\end{remark}

Observe that our
Definition of cointegral only makes sense in
the case where $H$ is finite
dimensional: indeed, we need a
dual basis of $H$ in order to make $H^*$
into a right
quasi-Hopf bimodule (see \ref{str2}). Also the equivalent
characterizations
from \prref{3.4} make no sense in the infinite
dimensional case,
as they involve the distinguished element, which can only
be
defined in the finite dimensional case. Nevertheless, the
cointegral has
a lot of applications in the finite dimensional
case (see
\cite{hn3}).\\
Following  \cite[Lemma 4.1]{bn2}, we now give an alternative
definition
for the space of coinvariants of a right quasi-Hopf $H$-bimodule
$M$:
\begin{equation}\label{scoinv}
M^{\ov {coH}}=\{n\in M\mid \rho
(n)=\sum x^1\cd n\cd
S(x^3_2X^3)f^1 \ot x^2X^1\b
S(x^3_1X^2)f^2\}.
\end{equation}
Using this definition, we will prove a
second
Structure Theorem for right quasi-Hopf bimodules. For $M\in
\hb $,
we define
\begin{equation}\label{overe}
\ove :\ M\ra
M,
\mbox{${\;\;\;}$}
\ove (m):=\sum m_{(0)}\cd \b
S(m_{(1)}),
\mbox{${\;\;\;}$$\forall m\in M$.}
\end{equation}
It follows
that
\begin{equation} \label{eove}
\ove (m)=\sum E(p^1\cd m)\cd
p^2,
\mbox{${\;\;}$}
E(m)=\sum X^1\cd \ove (m)\cd S(X^2)\a
X^3,
\mbox{${\;\;}$$\forall m\in M$.}%
\end{equation}
By \cite[Proposition
3.4]{hn3}, we have
\begin{equation}\label{fe}
E(m\cd h)=\va
(h)E(m),
\mbox{${\;\;}$}
E(h\cd E(m))=E(h\cd m),
\mbox{${\;\;}$}
h\cd
E(m)=\sum E(h_1\cd m)\cd h_2
\end{equation}
for all $m\in M$ and $h\in
H$.
and therefore the maps
\begin{equation}\label{in}
\ove:\ M^{{\rm
co}H}\to \ovc~~{\rm and}~~
E:\ \ovc\to M^{{\rm co}H}
\end{equation}
are
each others inverses. In the case of a classical Hopf algebra,
the maps $E$
and $\ove$ are equal to the identity on
$M^{{\rm co}H}=\ovc$.
Moreover, in
this case $M^{{\rm co}H}$ is
invariant under the left adjoint $H$-action
$h\tr m:=\sum h_1\cd
m\cd S(h_2)$, in the sense that $E(h\tr m)=h\tr E(m)$,
$h\in H$,
$m\in M$. In the quasi-Hopf case the projection $E$
generalizes
this property. More precisely, if we define $h\neg
m=E(h\cd m)$, then by
\cite[Proposition 3.4]{hn3}, we have that
$h\cd E(m)=\sum [h_1\neg E(m)]\cd
h_2$. Now, as in the Hopf algebra case,
we will prove that $\ovc $ is
invariant under the left
adjoint
$H$-action.

\begin{lemma}\lelabel{3.5}
Let $H$ be a quasi-Hopf
algebra, $M$ a right quasi-Hopf
$H$-bimodule
\begin{enumerate}
\item[(a)]
${\rm Im}(\ove )\subseteq \ovc $. If $n\in M$
then $n\in \ovc
\Leftrightarrow \ove (n)=n$.
\item[(b)] $\ovc $ is a left $H$-submodule of
$M$, where $M$ is
considered a left $H$-module via the left adjoint action,
that is
$h\tr m=\sum h_1\cd m\cd S(h_2)$, for all $h\in H$ and $m\in
M$.
\end{enumerate}
\end{lemma}

\begin{proof}
(a) Let $m\in M$ and $\d
=\sum \d ^1\ot \d ^2$ be given
by (\ref{gd}).
Then
\begin{eqnarray*}
\hspace*{-2cm}\rho (\ove (m))
&=&\sum m_{(0, 0)}\cd
\b _1S(m_{(1)})_1\ot
m_{(0, 1)}\b
_2S(m_{(1)})_2\\
{\rm(\ref{gdf},\ref{ca})}~~~~&=&\sum m_{(0, 0)}\cd \d
^1S(m_{(1)_2})f^1\ot
m_{(0, 1)}\d
^2S(m_{(1)_1})f^2
\\
{\rm(\ref{gd},\ref{q3},\ref{q5})}~~~~&=&\sum m_{(0, 0)}\cd
x^1\b
S((m_{(1)}x^3)_2X^3)f^1\ot  m_{(0, 1)}x^2X^1\b
S((m_{(1)}x^3)_1X^2)f^2\\
{\rm(\ref{hb2}, \ref{q1},\ref{q5})}~~~~&=&\sum
x^1\cd m_{(0)}\cd \b
S(x^3_2X^3m_{(1)})f^1\ot
x^2X^1\b
S(x^3_1X^2)f^2\\
&=&\sum x^1\cd \ove (m)\cd S(x^3_2X^3)f^1\ot
x^2X^1\b
S(x^3_1X^2)f^2.
\end{eqnarray*}
Therefore ${\rm Im} (\ove
)\subseteq \ovc $, implying immediately
one implication of the second
part.
Conversely, if $n\in \ovc $ then
\begin{eqnarray*}
\ove (n)
&=&\sum
n_{(0)}\cd \b S(n_{(1)})\\
{\rm(\ref{scoinv}, \ref{fgab})}~~~~&=&\sum
x^1\cd 
n\cd S(x^2X^1\b S(x^3_1X^2)\a x^3_2X^3)\\
{\rm(\ref{q5},
\ref{q6})}~~~~&=&n.
\end{eqnarray*}

(b) If $M$ is an $H$-bimodule then is
not hard to see that $M$ is
a left $H$-module via the $H$-action $h\tr
m=\sum h_1\cd m\cd
S(h_2)$, $h\in H$, $m\in M$. Using (\ref{eove}),
(\ref{fe}) and (\ref{qr1}) 
we have that $\ove (h\tr m)=h\tr \ove (m)$ for
all $h\in H$, $m\in M$. 
By part (a), we see  that $\ovc $ is a left
$H$-module under 
the left adjoint action.
\end{proof}

Let $H$ be a
quasi-Hopf algebra and $M$ a right quasi-Hopf
$H$-bimodule. We already have
seen that the map
$\ove:\ M^{{\rm co}H}\to\ovc $ is an isomorphism. It is
also
$H$-linear since
\begin{eqnarray*}
\ove (h\neg n)
&=&\ove (E(h\cd
n))\\
{\rm(\ref{eove}, \ref{fe})}~~~~&=&\sum E(p^1h\cd n)\cd
p^2\\
{\rm(\ref{qr1}, \ref{fe})}~~~~&=&\sum h_1\cd E(p^1\cd n)\cd
p^2S(h_2)\\
{\rm(\ref{eove})}~~~~&=&h\tr \ove (n)
\end{eqnarray*}
for all
$h\in H$ and $n\in M^{{\rm co}H}$.\\
$\ovc $ is a left $H$-module, so, by
\cite[Lemma 3.2]{hn3},
$\ovc \ot H$ becomes a right quasi-Hopf $H$-bimodule
with the
following structure:

\begin{equation}\label{ehb} a\cd (n\ot h)\cd
b=\sum a_1\tr n\ot a_2hb,
\mbox{${\;\;\;}$} \rho '(n\ot h)=\sum x^1\tr n\ot
x^2h_1\ot
x^3h_2
\end{equation}
for all $a, b, h\in H$ and $n\in \ovc $. We
can now state the second
Structure Theorem for right quasi-Hopf
bimodules.

\begin{theorem}\thlabel{3.6}
Let $H$ be a quasi-Hopf algebra
and $M$ a right quasi-Hopf
$H$-bimodule. Consider $\ovc \ot H$ as a right
quasi-Hopf
$H$-bimodule as in (\ref{ehb}). The
map
\begin{equation}\label{onu}
\ov {\nu }: \ovc \ot H\ra
M,
\mbox{${\;\;}$}
\ov {\nu }(n\ot h)=\sum X^1\cd n\cd S(X^2)\a
X^3h
\end{equation}
is an isomorphism of quasi-Hopf $H$-bimodules. The
inverse of
$\ov {\nu }$ is given by the
formula
\begin{equation}\label{ionu}
\ov {\nu }^{-1}(m)=\sum \ove
(m_{(0)})\ot m_{(1)}.
\end{equation}
\end{theorem}

\begin{proof}
We have
seen that $\ovc\cong M^{{\rm co}H}$ are isomorphic as
left $H$-modules, and
therefore $\ovc\ot H\cong M^{{\rm co}H}\ot H$
as quasi-Hopf $H$-bimodules
(in both cases, the structure is
determined by \cite[Lemma 3.2]{hn3}). From
the Hausser-Nill
Structure Theorem (\cite[Theorem 3.8]{hn3}),
it follows
that $M\cong M^{{\rm co}H}\ot H$
as quasi-Hopf $H$-bimodules. Thus we find
that $\ovc\ot H\cong M$
as quasi-Hopf $H$-bimodules, and it is
straightforward to verify
that the connecting isomorphism is exactly
$\ov{\nu}$.
\end{proof}

\begin{remark}\rm \thref{3.6} can be proved also
in a direct way,
without using the Hausser-Nill Structure Theorem. This
proof
is straightforward, but long and technical. However, it has
the
advantage that it is independent of the bijectivity of the
antipode,
which makes it more general in the infinite dimensional case.
At
first sight,
the Hausser-Nill definition of coinvariants makes use of the
bijectivity
of the antipode; however, a careful inspection of their
definition and
their proof of the Structure Theorem shows that we do not
really need the
bijectivity.
\end{remark}

The antipode of a finite
dimensional quasi-Hopf algebra is bijective,
so its dual $H^*$ is a right
quasi-Hopf $H$-bimodule with
structures defined in (\ref{str1}),
(\ref{str2}). The coinvariants
$\ov{\lambda} \in H^{*\ov {coH}}$ are called
left alternative cointegrals
on $H$,
and the space of left alternative
cointegrals is denoted by $\ov{\cal L}=
H^{*\ov
{{\rm co}H}}$. From
\thref{3.6}, we obtain immediately the following
result:

\begin{theorem}\thlabel{3.7}
Let $H$ be a finite dimensional
quasi-Hopf algebra. Then
$dim_k\ov{\cal L}=1$.
\end{theorem}

Applying
(\ref{q3}) and (\ref{q5}),
we find that $\ov{\lambda} \in H^*$ is an
alternative left cointegral
if and only if
\begin{equation}\label{Int}
\sum
\ov{\lambda} (V^2h_2U^2)V^1h_1U^1=\sum
\ov{\lambda} (\smi
(X^1_1p^1)hS(S(X^3)f^1))X^1_2p^2S(X^2)f^2
\end{equation}
for all $h\in H$.
(\ref{Int}) can be used to extend the definition of
left alternative
cointegral to infinite dimensional quasi-Hopf algebras.

\section{Integrals
for dual quasi-Hopf
algebras}\selabel{4}
\setcounter{equation}{0}
Sullivan's Theorem \cite{su}
asserts that the space of left integrals
on a Hopf algebra has dimension at
most one. Various
new proofs of this result have been given in recent
years, see \cite{st},
\cite{bdgn}, \cite{dnt} and \cite{ra}. The aim of
this Section is to
give a proof of the uniqueness of integrals on a dual
quasi-Hopf
algebra. Our approach is based on the methods developed in
\cite{bdgn}.\\
Throughout, $A$ will be a dual quasi-bialgebra or a dual
quasi-Hopf
algebra. Following \cite{maj}, a dual quasi-bialgebra
$A$ is a
coassociative coalgebra $A$ with comultiplication
$\Delta $ and counit $\va
$ together with coalgebra morphisms
$M:\ A\ot A\ra A$ (the multiplication;
we write $M(a\ot b)=ab$)
and $u:\ k\ra A$ (the unit; we write $u(1)=1$),
and an
invertible element $\v \in (A\ot A\ot A)^*$ (the reassociator),
such
that for all $a, b, c, d\in A$ the following relations
hold:
\begin{equation}\label{dq1}
\sum a_1(b_1c_1)\v (a_2, b_2, c_2)=
\sum
\v (a_1, b_1,
c_1)(a_2b_2)c_2,
\end{equation}
\begin{equation}\label{dq2}
1a=a1=a,
\end{equation}
\begin{equation}\label{dq3}
\sum \v (a_1, b_1, c_1d_1)\v (a_2b_2, c_2,
d_2)=
\sum \v (b_1, c_1, d_1)\v (a_1, b_2c_2, d_2)
\v (a_2, b_3, c_3),
\end{equation}
\begin{equation}\label{dq4}
\v (a, 1, b)=\va (a)\va (b).
\end{equation}
$A$ is called a dual quasi-Hopf algebra if, moreover,
there exist an
antimorphism $S$ of the coalgebra $A$ and elements $\a , \b
\in H^*$
such that, for all $a\in A$:
\begin{equation}\label{dq5}
\sum
S(a_1)\a (a_2)a_3=\a (a)1,
\mbox{${\;\;\;}$}
\sum a_1\b (a_2)S(a_3)=\b
(a)1,
\end{equation}
\begin{equation}\label{dq6}
\sum \v (a_1\b (a_2),
S(a_3), \a (a_4)a_5)=\sum
\v ^{-1}(S(a_1), \a (a_2)a_3, \b (a_4)S(a_5))=\va
(a).
\end{equation}
It follows from the axioms that $\a (1)\b (1)=1$,
so we
can assume that $\a (1)=\b (1)=1$ and,
$S(1)=1$. Moreover (\ref{dq3}) and
(\ref{dq4}) imply
\begin{equation}\label{dq7}
\v (1, a, b)=\v (a, b, 1)=\va
(a)\va (b),
\mbox{${\;\;}$$\forall a, b\in A$.}
\end{equation}
If $A=(A, M,
u, \v , S, \a , \b )$ is a dual quasi-Hopf algebra,
then $A^{\rm op, cop}$
is also a dual quasi-Hopf algebra. The structure
maps are$$
\v _{\rm op,
cop}(a, b, c)=\v (c, b, a),\; S_{\rm op, cop}=S,\; \a _{\rm
op, cop}
=\b \;
\mbox{and}\; \b _{\rm op, cop}=\a . $$
If $A$ is
a dual quasi-bialgebra,
then $A^*$ is an algebra, with multiplication
given by convolution, and
unit $\va$.

\begin{definition}\delabel{4.1}
A map $T\in A^*$ is called a
left integral on
the dual quasi-bialgebra $A$ if $a^*T=a^*(1)T$ for
any
$a^*\in A^*$. Left integrals on $A^{\rm op, cop}$ are called
right
integrals on $A$.
\end{definition}

The set of left (right) integrals on
$A$ is denoted
by $\int _l$ ($\int _r$).
We keep the same notation as in
\seref{3}
but we will specify every time which kind of
integral we are
using.
It is clear that $T\in A^*$ is a left integral if and only if
$\sum
T(a_2)a_1=T(a)1$ for all $a\in A$. $\int _l$ is
a subspace of $A^*$.
Moreover, as in the Hopf algebra case, $\int _l$ is
a two-sided ideal in the algebra $A^*$.\\
Let $A^{*rat}$ be the left rational
part of $A^*$. $A^{*rat}$ is the sum of
rational left ideals of the algebra $A^*$, see \cite[Ch. II]{sw}.
Note that $A^{*rat}\subseteq A^*$ and
\begin{eqnarray}
&&\hspace*{-15mm}a^*\in A^{*rat}\Longleftrightarrow\nonumber\\
&&\exists \; (a_i^*)_{i=\ov {1, n}}\subseteq  A^*
~~{\rm and}~~
(a_i)_{i=\ov {1, n}}\subseteq A
~{\rm such~that}~
b^*a^*=\sum \limits _{i=1}^nb^*(a_i)a_i^*,
~~\forall b^*\in A^*.\label{rat}
\end{eqnarray}
It follows that $\int _l\subseteq A^{*rat}$. In particular,
if  $A^{*rat}=0$, then
$\int _l=0$.\\
Later in this Section we will show that
the left and right rational parts of $A^*$
are equal, justifying our notation.
As in the Hopf algebra case, we first describe the
connection between $A^{*rat}$ and $\int _l$. It is well-known
that $A^{*rat}$ is a rational left $A^*$-module, and this induces
a right $A$-comodule structure on $A^{*rat}$ defined by
\begin{equation}\label{rat2}
\rho :\ A^{*rat}\ra A^{*rat}\ot A,
\mbox{${\;}$}
\rho (a^*)=\sum a^*_{(0)}\ot a^*_{(1)}\Leftrightarrow
b^*a^*=\sum b^*(a^*_{(1)})a^*_{(0)},
\end{equation}
for all $b^*\in A^*$. This can be rewritten as follows
\begin{equation}\label{rat3}
\rho (a^*)=\sum a^*_{(0)}\ot a^*_{(1)}\Leftrightarrow
\sum a^*(a_2)a_1=\sum a^*_{(0)}(a)a^*_{(1)},\;
\forall a\in A.
\end{equation}
Now, define the map $\sigma :\ A\ot A\ra A^*$ by
\begin{equation}\label{st}
\sigma (a\ot b)(c)=\v (c, a, b),
\mbox{${\;}$$\forall a, b, c\in A$.}
\end{equation}
$\sigma $ is
convolution invertible, the
inverse $\sigma ^{-1}$ is given by
$\sigma
^{-1}(a\ot b)(c)=\v ^{-1}(c, a, b)$, for all $a, b, c\in A$.
We introduce
the following notation, for $a^*\in A^*$ and $a\in
A$:
$$
a^*\leftharpoondown a=S(a)\rh a^*\in A^*.
$$
For all $a^*\in
A^{*rat}$, we define $P^*(a^*)$ by
\begin{equation}\label{ps}
P^*(a^*)=\sum
\b (a^*_{(1)})\a (S(a^*_{(3)}))
\sigma (S^2(a^*_{(4)})\ot
S(a^*_{(2)}))(a^*_{(0)}\leftharpoondown
S(a^*_{(5)})).
\end{equation}
We
now claim that
$P^*(a^*)\in \int _l$,
for all $a^*\in A^{*rat}$. Indeed,
for all $b^*\in A^*$, we
calculate
\begin{eqnarray*}
&&\hspace*{-3cm}\langle b^*(1)P^*(a^*),
a\rangle \\
{\rm(\ref{ps}, \ref{dq5})}~~~~&=&
\sum \langle b^*, a^*_{(1)}\b
(a^*_{(2)})S(a^*_{(3)})\rangle \a
(S(a^*_{(5)}))\v (a_1, S^2(a^*_{(6)}),
S(a^*_{(4)}))\\
&&~~~~\langle a^*_{(0)}, a_2S^2(a^*_{(7)})\rangle
\\
{\rm(\ref{rat3}, \ref{dq1})}~~~~&=&\sum \b (a^*_{(1)})
\langle b^*,
a_1[S^2(a^*_{(5)})\a (S(a^*_{(4)}))S(a^*_{(3)})]\rangle
\v (a_2,
S^2(a^*_{(6)}), S(a^*_{(2)}))\\
&&~~~~\langle a^*_{(0)},
a_3S^2(a^*_{(7)})\rangle \\
{\rm(\ref{dq5}, \ref{ps})}~~~~&=&\sum
b^*(a_1)\langle P^*(a^*), a_2\rangle 
=\langle b^*P^*(a^*), a\rangle .
\end{eqnarray*}

\begin{proposition}\prlabel{4.2}
Let $A$ be a dual
quasi-Hopf algebra and
$\sigma :\ A\ot A\ra A^*$ the map defined in
(\ref{st}).
Then $\theta ^*:\ \int _l\ot A\ra A^{*rat}$, given
by
\begin{equation}\label{ts}
\theta ^*(T\ot a)=\sum \sigma (S(a_5), \a (a_6)a_7)
(T\leftharpoondown a_4)
\sigma ^{-1}(S(a_3)\ot \b
(S(a_2))S^2(a_1)),
\end{equation}
is an isomorphism of right
$A$-comodules.
\end{proposition}

\begin{proof}
In the situation where $A$
is finite dimensional, the proof follows from
\thref{2.2} by duality. This
is why we restrict to proving that
$\theta ^*$ is well-defined and has an
inverse, leaving other details
to the reader. For $a^*\in A^*$, $T\in \int
_l$ and $a\in A$, we compute that
\begin{eqnarray*}
&&\hspace*{-2cm}\langle
a^*\theta ^*(T\ot a), b\rangle \\
&=&\sum a^*(b_1)\v (b_2, S(a_5), \a
(a_6)a_7)T(b_3S(a_4))
\v ^{-1}(b_4, S(a_3), \b
(S(a_2))S^2(a_1))\\
{\rm(\ref{dq5})}~~~~&=&\sum \langle a^*, b_1(S(a_6)\a
(a_7)a_8)\rangle \v
(b_2, S(a_5), a_9)T(b_3S(a_4))\\
&&~~~~\v ^{-1}(b_4,
S(a_3), \b (S(a_2))S^2(a_1))\\
{\rm(\ref{dq1})}~~~~&=&\sum \langle (a_8\rh
a^*)T, b_2S(a_4)\rangle 
\v (b_1, S(a_5), \a (a_6)a_7)\\
&&\v ^{-1}(b_3,
S(a_3), \b (S(a_2))S^2(a_1))\\
&=&\sum a^*(a_2)\langle \theta ^*(T\ot a_1), 
b\rangle 
\end{eqnarray*}
hence
\begin{equation}\label{F}
a^*\theta ^*(T\ot a)=\sum a^*(a_2)\theta ^*(T\ot a_1).
\end{equation}
From (\ref{rat}) it is follows that $\theta ^*(T\ot a)\in A^{*rat}$,
as needed.\\
$\int _l\ot A$ is a right $A$-comodule with structure induced by
the comultiplication on $A$, and it follows then from (\ref{rat2})
that $\theta ^*$ is right $A$-colinear.\\
We claim that the inverse of $\theta^*$ is given by
\begin{equation}\label{ints}
\theta ^{*-1}:\ A^{*rat}\ra \int _l\ot A,
\mbox{${\;}$}
\theta ^{*-1}(a^*)=\sum P^*(a^*_{(0)})\ot a^*_{(1)},
\end{equation}
where $\rho (a^*)=\sum a^*_{(0)}\ot a^*_{(1)}\in A^{*rat}\ot A$
is defined as in (\ref{rat2}).
It is clear that $\theta ^{*-1}$ is well defined. To show
that $\theta ^*$ and $\theta ^{*-1}$ are each others inverses, we need the
following equalities, for any $T\in \int _l$ and $a, b\in A$:
\begin{eqnarray}
&&\hspace*{-2cm}\sum \v (a_2, S(b_2), \a (b_3)b_4)T(a_3S(b_1))a_1\nonumber\\
&=&\sum \v (a_1, S(b_2), \a (b_3)b_4)T(a_2S(b_1))b_5,\label{F1}\\
&&\hspace*{-2cm}
T(aS(b))=\sum \v (a_1, S(b_2)\a (b_3), b_4\b (b_5))T(a_2S(b_1))\nonumber\\
&=&\sum \v (\b (a_1)a_2, S(b_2), \a (b_3)b_4)T(a_3S(b_1))\label{F2}.
\end{eqnarray}
Note that these formulas are the formal
duals of (\ref{f1}) and (\ref{f2}).
\end{proof}

\begin{corollary}\colabel{4.3}
Let $A$ be a dual quasi-Hopf algebra. Then $A^{*rat}=0$ if and only if
$\int _l=0$.
\end{corollary}

\begin{corollary}\colabel{4.4}
Assume that a dual quasi-Hopf algebra $A$ has a non-zero left integral.
Then the the antipode $S$ is injective. If $A$ is finite dimensional,
then the left integral space has dimension 1, and the antipode $S$ is
bijective.
\end{corollary}

\begin{proof}
Let $T$ be a non-zero left integral, and assume
that
$S(a)=0$. Let $\ov {S}$ be the algebra
antimorphism dual to $S$. If
$\theta ^*$ is the map defined
by (\ref{ts}) then for all $a^*\in A^*$ and
$b\in A$ we have:
\begin{eqnarray*}
&&\hspace*{-2cm}\langle \theta ^*(T\ot
\b \ov {S}(a^*)\rh a), b\rangle =
\sum \b (a_2)a^*(S(a_3))\langle \theta
^*(T\ot a_1), b\rangle \\
&=&\sum \b (a_8)a^*(S(a_9))\v (b_1, S(a_5), \a
(a_6)a_7)T(b_2S(a_4))\\
&&\v ^{-1}(b_3, S(a_3), \b (S(a_2))S^2(a_1))\\
{\rm
(\ref{F2})}~~~~&=&\sum a^*(S(a_5))T(b_1S(a_4))\v ^{-1}(b_2, S(a_3),
\b
(S(a_2))S^2(a_1))=0.
\end{eqnarray*}
Since $\theta ^*$ is bijective and
$T\not=0$ it is follows that
$\b \ov {S}(a^*)\rh a=0$, for any $a^*\in
A^*$.
Thus $\sum \b (a_2)a_1\ot S(a_3)=0$ and therefore
$$
\sum \b
(a_5)a_1\ot a_2\ot a_3\ot a_4\ot S(a_6)=0.
$$
By (\ref{dq6}) we obtain that
$a=\sum
\v ^{-1}(S(a_2), \a (a_3)a_4, \b (a_5)S(a_6))a_1=0$,
and $S$ is
injective.\\
If $A$ is finite dimensional then $A^{*rat}=A^*$. We obtain
that
$\theta ^*:\ \int _l\ot A\ra A^*$ is an isomorphism of right
$A$-comodules.
The final assertion is then obvious.
\end{proof}

Let $C$ be
a coalgebra and $C^*$ the dual algebra. Then
$C$ is a left (right)
$C^*$-module under the left (right)
action $\rh $ ($\lh $) of $C^*$ on $C$
given by $c^*\rh c=\sum
c^*(c_2)c_1$ ($c\lh c^*=\sum c^*(c_1)c_2$), for all
$c^*\in C^*$
and $c\in C$. Recall that $C$ is called a left
(right) quasi-co-Frobenius coalgebra 
(shortly QcF coalgebra) if there 
exists an injective morphism of left 
(right) $C^*$-modules from $C$
to a free left (right) $C^*$-module. The coalgebra 
$C$ is called left (right) co-Frobenius if there 
exists a monomorphism of (left) right 
$C^*$-modules from $C$ to $C^*$ or, equivalently, 
if there exists a 
bilinear form ${\bf b}:\ C\ot C\ra k$ which is left 
(right) non-degenerated 
and $C^*$-balanced, i.e. if ${\bf b}(c, x)=0$ for 
any $c\in C$ (resp. ${\bf b}(x, c)=0$ for any
$c\in C$) then $x=0$ and 
${\bf b}(x\lh c^*, y)={\bf b}(x, c^*\rh y)$ 
for any $x, y\in C$, $c^*\in C^*$.\\ 
Finally, $C^*$ is a 
left (right) $C^*$-module, so we can consider the 
left (right) rational 
part of $C^*$. We will denote this 
$C^*$-submodule of $C^*$ by $C^{*rat}_l$ 
(resp. $C^{*rat}_r$). 
$C$ is called a left (right) semiperfect 
coalgebra if 
the category $^C{\cal M}$ (${\cal M}^C$) of left 
(right) $C$-comodules has 
enough projectives or, equivalently, 
$C^{*rat}_r$ ($C^{*rat}_l$) is a dense 
subset of $C^*$ (in  the finite 
topology, see \cite{ab}, \cite{dnr} for the 
Definition). Following 
\cite{bdgn}, if $C$ is a left and right semiperfect 
coalgebra then 
$C^{*rat}_l=C^{*rat}_r:=C^{*rat}$, $C^{*rat}$ is dense in
$C^*$ 
and $C$ is projective generator in the 
categories $^C{\cal M}$  and ${\cal M}^C$.\\
It is known that a left (right)
co-Frobenius coalgebra
is a left (right) QcF coalgebra and a left
(right) QcF coalgebra is a left
(right) semiperfect coalgebra, but
the converse implications are not true.
There are examples showing
that none of the three concepts is left-right
symmetric. However, if
$H$ is a Hopf algebra all these concepts are
equivalent to $H$
having non-zero left (or right) integrals. We prove now
that a
similar result holds for dual quasi-Hopf algebras. Except for
the
implication $(i)\Rightarrow (ii)$, the proof is identical to the
proof
of \cite[Theorem 5.3.2]{dnr}, so we omit it
here.

\begin{theorem}\thlabel{4.5}
Let $A$ be a dual quasi-Hopf algebra.
Then the following
assertions are equivalent:\\
\hspace*{1cm}(i) $A$ has a
non-zero integral.\\
\hspace*{1cm}(ii) $A$ is a left co-Frobenius
coalgebra.\\
\hspace*{1cm}(iii) $A$ is a left QcF
coalgebra.\\
\hspace*{1cm}(iv) $A$ is a left semiperfect
coalgebra.\\
\hspace*{1cm}(v) $A$ has a non-zero right
integral.\\
\hspace*{1cm}(vi) $A$ is a right co-Frobenius
coalgebra.\\
\hspace*{1cm}(vii) $A$ is a right QcF
coalgebra.\\
\hspace*{1cm}(viii) $A$ is a right semiperfect
coalgebra.\\
\hspace*{1cm}(ix) $A$ is a generator in the category $^A{\cal
M}$ (or
in ${\cal M}^A$).\\
\hspace*{1cm}(x) $A$ is a projective object in
the category $^A{\cal M}$
(or in ${\cal
M}^A$.)
\end{theorem}

\begin{proof}$\un{(i)\Rightarrow (ii)}$ Let $T\in
A^*$ be a non-zero
left integral. We define a bilinear form ${\bf b}: A\ot
A\ra k$ as
follows:
\begin{equation}\label{X}
{\bf b}(a, b)=\theta ^*(T\ot
b)(a),
\mbox{${\;}$$\forall a, b\in A$,}
\end{equation}
where $\theta ^*$
is the map defined as in (\ref{ts}). Then, for
all $a, b\in A$ and $a^*\in
A^*$ we compute:
\begin{eqnarray*}
{\bf b}(a\lh a^*, b)
&=&\sum
a^*(a_1){\bf b}(a_2, b)\\
{\rm (\ref{X})}~~~~&=&\langle a^*\theta ^*(T\ot
b), a\rangle \\
{\rm (\ref{F}, \ref{X})}~~~~&=&\sum a^*(b_2){\bf b}(a,
b_1)
={\bf b}(a, a^*\rh b)
\end{eqnarray*}
proving that ${\bf b}$ is
$C^*$-balanced. Now we prove that ${\bf b}$
is left non-degenerate. If
$x\in A$ such that ${\bf b}(a, x)=0$ for
any $a\in A$ then $\theta ^*(T\ot
x)=0$, and it follows from \prref{4.2}
that $x=0$. Thus ${\bf b}$ is
$C^*$-balanced and left non-degenerate,
so $A$ is a left co-Frobenius
coalgebra.
\end{proof}

By \cite[Proposition 2.2]{dnt}, any subcoalgebra
$D$ of
a left semiperfect coalgebra $C$ is itself left semiperfect.
As a
consequence of \thref{4.5}, we therefore
have

\begin{corollary}\colabel{4.6}
Let $A$ be a dual quasi-Hopf algebra
with non-zero integrals. Then any
dual quasi-Hopf subalgebra $B$ of $A$
(i.e. a subcoalgebra $B$ of $A$
which is closed under multiplication of
$A$, $1\in B$ and
$S(B)\subseteq B$) has non-zero
integrals.
\end{corollary}

We proceed with a proof of the uniqueness
of
integrals for dual quasi-Hopf
algebras. First, by \cite[Remark
5.4.3]{dnr}, if $C$ is a left
and right co-Frobenius coalgebra, and $M$ is
a finite dimensional
right $C$-comodule, then ${\rm dim}_k{\rm
Hom}_{C^*}(C, M)\leq {\rm dim}_kM$
(Recall that $M$ and $N$ are in a
natural way left
$C^*$-modules).

\begin{proposition}\prlabel{4.7}
Let $A$
be a dual quasi-Hopf algebra with non-zero
integral. Then ${\rm dim}_k\int
_l={\rm dim}_k\int _r=1$.
\end{proposition}

\begin{proof}
$A$ is a left
and right co-Frobenius dual quasi-Hopf algebra
by \thref{4.5}. $k$ is a
right $A$-comodule, and
$\int _r={\rm Hom}_{A^*}(A, k)$. Thus, $\int _r$
has dimension at most $1$, so
this has dimension precisely $1$ since $\int
_r\not=0$. Now, if
we replace the dual quasi-Hopf algebra $A$ with $A^{\rm
op, cop}$ then we
obtain that ${\rm dim}_k\int _l=1$.
\end{proof}

A result
of Radford \cite{re} asserts that a co-Frobenius Hopf
algebra has bijective
antipode. It is not known if a
similar result holds for a co-Frobenius dual
quasi-Hopf algebra.
To this end, we would need the dual version of
\leref{2.3},
without the assumption that the antipode is
bijective.
However, we know that the antipode of
a co-Frobenius dual
quasi-Hopf algebra is injective (\coref{4.4}). We apply
this to prove the
following result:

\begin{lemma}\lelabel{4.8}
Let $A$ be a dual quasi-Hopf
algebra with a non-zero
left integral $T$ and
antipode $S$. Then $T\circ S$
is a non-zero right integral on $A$.
In particular, $\ov {S}(\int _l)=\int
_r$.
\end{lemma}

\begin{proof} We omit the proof, since it is identical to
proof of
\cite[Lemma 5.4.4]{dnr}.
\end{proof}

The right integral $T\circ
S$ can be described more explicitly.
Indeed, let $A$ be a dual-quasi
Hopf algebra with non-zero left integral $T$.
Then, for all $a^*\in A^*$,
$Ta^*$ is also a left integral and
$Ta^*=\chi (a^*)T$ for some map $\chi $
from $A^*$ to $k$. Since
$Ta^*b^*=\chi (a^*b^*)T=\chi (a^*)\chi (b^*)T$ and
$T=\chi (\va )T$,
$\chi $ is an algebra map.\\
Let $G(C)$ be the set of
grouplike elements of a coalgebra $C$.
From (\ref{dq6}) and (\ref{dq5}), it
follows that a grouplike
element $g$ in a dual quasi-Hopf algebra is
invertible,
and $g^{-1}=S(g)$.\\
The following result generalizes
\cite[Proposition 1.3]{bdr}.

\begin{proposition}\prlabel{4.9}
Let $A$ be a
dual quasi-Hopf algebra with a non-zero left integral $T$.
Then there
exists a grouplike element $g$ in $A$ such that\\[1mm]
${\;\;\;}$(i)
$Ta^*=a^*(g)T$, for any $a^*\in A^*$;\\[1mm]
${\;\;\;}$(ii) $\ov
{S}(T)=\Lambda (g\rh T)$, where $\Lambda \in A^*$ is given by 
$\Lambda
(b)=\sum \a (b_2)\v ^{-1}(S(b_1), b_3, g)$, for all $b\in
A$.
\end{proposition}

\begin{proof}
Since $T\not=0$, there exists $a\in A$
such that
$T(a)=1$. We will show that
$$g=a\lh T=\sum T(a_1)a_2$$
has the
required properties. As in \seref{3}, we
denote by $\rh $ and $\lh $ the
usual left and right action respectively
of $A$ on $A^*$, that
is
$$
\langle b\rh a^*, c\rangle =\langle a^*,
cb\rangle
\mbox{${\;\;}$and${\;\;}$}
\langle a^*\lh b, c\rangle =\langle
a^*, bc\rangle .
$$
It is easy to show that
\begin{eqnarray}
(a^*b^*)\lh
b&=&\sum (a^*\lh b_1)(b^*\lh b_2);\label{ma}\\
b\rh (a^*b^*)&=&\sum (b_1\rh
a^*)(b_2\rh b^*).\label{ma2}
\end{eqnarray}
Let $\chi\in A^{**}$ be defined
as
above. For all $a^*\in A^*$, we have that
$$
\chi (a^*)=\chi
(a^*)T(a)=\langle Ta^*, a\rangle =\sum T(a_1)a^*(a_2)=
\langle a^*, \sum
T(a_1)a_2\rangle =a^*(g),
$$
and therefore $Ta^*=\chi (a^*)T=a^*(g)T$, for any $a^*\in A^*$.
$g$ is a grouplike element, because $\chi $ is an algebra map.\\
(ii) Define $\ov {\Lambda }: A\ra A^*$ given by
$\ov {\Lambda }(b)(c)=\sum \a (b_2)\v ^{-1}(S(b_1), b_3, c)$ for any
$b, c\in A$. Is follows that $\Lambda (b)=\ov {\Lambda }(b)(g)$, for any
$b\in A$,
and $\ov {\Lambda }(1)=\va $. Moreover, by (\ref{dq1}), (\ref{dq5}) and (i)
we have that
\begin{equation}\label{Y}
\sum \ov {\Lambda }(b_2)[(T\lh S(b_1))\lh b_3]=T\ov {\Lambda }(b)=\Lambda (b)T
\mbox{${\;\;\;}$$\forall $ $b\in B$.}
\end{equation}
Now, the relation in (ii) follows from the following computation,
for all $b\in A$:
\begin{eqnarray*}
\hspace*{-2cm}\langle \Lambda (g\rh T), b\rangle &=&\sum \Lambda
(b_1)\langle T\lh
b_2, g\rangle T(a)\\
{\rm (i)}~~~~&=&\sum \Lambda (b_1)\langle T(T\lh b_2),a\rangle \\
{\rm (\ref{Y}, \ref{ma})}~~~~&=&
\sum \langle \ov {\Lambda }(b_2)\{[(T\lh S(b_1))T]\lh b_3\},
a\rangle \\
&=&\sum T(S(b_1))\langle \ov {\Lambda }(b_2)(T\lh b_3), a\rangle
\mbox{${\;\;\;}$since $T\in \int _l$}\\%
&=&\langle \ov {\Lambda }(1)T, a\rangle \langle \ov {S}(T), b\rangle
\mbox{${\;\;\;}$ since $\ov {S}(T)\in \int _r$, see \leref{4.8}}\\
&=&\langle \ov {S}(T), b\rangle .
\end{eqnarray*}
\end{proof}
${\;\;\;}$

We conclude this paper with a Maschke type theorem for dual
quasi-Hopf algebras. Let $C$ be a coalgebra and $M$ a right
$C$-comodule. Recall that $M$ is called cosemisimple if for any
subcomodule $N$ of $M$ there exists a $C$-colinear map $\pi :\ M\ra
N$ such that $\pi \circ i=id$, where $i:\ N\ra M$ is the inclusion
map. $C$ is called cosemisimple if any right $C$-comodule is
cosemisimple.

\begin{theorem}\thlabel{4.10}
For a dual quasi-Hopf algebra $A$,the following
statements are equivalent:
\begin{enumerate}
\item[(i)] $A$ is a cosemisimple coalgebra;
\item[(ii)] $A$, viewed as a right $A$-comodule via comultiplication,
is cosemisimple as an $A$-comodule;
\item[(iii)] there exists a right integral $T\in A^*$ such that
$T(1)=1$;
\item[(iv)] there exists a left integral $T\in A^*$ such that $T(1)=1$.
\end{enumerate}
\end{theorem}

\begin{proof}
$\un{(i)\Rightarrow (ii)}$. is trivial.\\
$\un{(ii)\Rightarrow (iii)}$.
Let $u:\ k\ra A$ be the unit map of $A$. $k$ is trivially a right
$A$-comodule. Since $u$ is injective there
exists an $A$-colinear map $T:\ A\ra k$ such that $T(1)=1$. This
means that $T$ is a right integral for $A$ and $T(1)=1$.\\
$\un{(iii)\Rightarrow (iv)}$.
Let $T$ be a right integral on $A$ such
that $T(1)=1$. Then $T$ is a left integral on $A^{\rm op, cop}$ and
$T\circ S$ is a right integral on $A^{\rm op, cop}$, by \leref{4.8}.
Therefore, $T\circ S$ is a left integral for $A$, and since
$T(1)=1$ it follows that $(T\circ S)(1)=1$.\\
$\un{(iv)\Rightarrow (i)}$.
 Let $T\in \int _l$ be such that $T(1)=1$ and
$M$ a right $A$-comodule.
Then $M\ot A$
is also a right $A$-comodule right $A$-action $id \ot \Delta $.
$\rho _M:\ M\to M\ot A$ is an injective right $A$-colinear map. Moreover,
the map $\omega _M:\ M\ot A\ra M$ given by
\begin{equation}\label{om}
\omega _M(m\ot a)=\sum \a (a_3)\v (m_{(1)}, S(a_2), a_4)
\b (m_{(3)})T(m_{(2)}S(a_1))m_{(0)},
\end{equation}
 is right $A$-colinear (by (\ref{F1})), and $\omega _M\circ \r_M=id$
(by (\ref{dq5}) and (\ref{dq6})).\\
Now, if $N$ is an $A$-subcomodule of $M$, then there is a
$k$-linear map $\tilde {\pi }:\ M\ra N$ such that $\tilde {\pi
}\circ i=id$, where $i$ is the inclusion map $i:\ N\ra M$. We
define $\pi :\ M\ra N$, $\pi =\omega _N\circ (\tilde {\pi }\ot
id)\circ \rho _M$, where $\omega _N$ is the corresponding map
(\ref{om}) for the right $A$-comodule $N$. Is not hard to see that
$\tilde {\pi }\ot id$ is $A$-colinear, hence $\pi $ is $A$-colinear.
It also follows from
the above considerations that $\pi \circ i=id$, proving that
$A$ is a cosemisimple coalgebra.
\end{proof}

\begin{center}
{\bf ACKNOWLEDGEMENTS}
\end{center}
The authors thank the referee for his helpful comments, which
improved a first version of this paper. They also would like to thank
Peter Schauenburg for telling them about his approach to quasi-bialgebras
during Hopf algebra meetings held in Brussels and Toronto during the spring and
fall of 2002.

\end{document}